\documentclass[10pt]{amsart}

\usepackage{amsmath}
\usepackage{amssymb}
\usepackage{amsthm}
\usepackage{amscd}

\newtheorem{theorem}{Theorem}
\newtheorem{lemma}[theorem]{Lemma}
\newtheorem{conjecture}[theorem]{Conjecture}
\newtheorem{proposition}[theorem]{Proposition}
\newtheorem{corollary}[theorem]{Corollary}
\numberwithin{theorem}{section}

\theoremstyle{definition}

\newtheorem{example}[theorem]{Example}
\newtheorem{remark}[theorem]{Remark}

\newtheorem*{acknowledgments}{Acknowledgments}

\numberwithin{equation}{section}

\newcommand{\Z}{{\mathbb Z}}
\newcommand{\Aff}{\hbox{Aff}}
\newcommand{\ot}{\otimes}
\newcommand{\geh}{\mathfrak{g}}

\newcommand\cd{\cdot}
\newcommand\la{\lambda}
\renewcommand\P{{\mathcal P}}

\begin{document}

\begin{center}
{\Large{\bf Crystal Interpretation of}}
\vspace{1mm}\\
{\Large{\bf Kerov-Kirillov-Reshetikhin Bijection}}\vspace{2mm}\\
{\large Atsuo Kuniba, Masato Okado, Reiho Sakamoto,}\vspace{1mm}\\
{\large Taichiro Takagi and Yasuhiko Yamada}\vspace{3mm}
\end{center}

\begin{quotation}
{\small
A}{\tiny BSTRACT:}
{\small
The Kerov-Kirillov-Reshetikhin (KKR) bijection is
the crux in proving fermionic formulas.
It is defined by a combinatorial algorithm on 
rigged configurations and highest paths.
We reformulate the KKR bijection as a 
vertex operator by purely using
combinatorial $R$ in crystal base theory.
The result is viewed as a nested Bethe ansatz at $q=0$ 
as well as the direct and 
the inverse scattering (Gel'fand-Levitan) map 
in the associated soliton cellular automaton.}
\end{quotation}


\section{Introduction}\label{sec:1}

Among many approaches to quantum integrable systems, 
Bethe ansatz \cite{Be} stands as a most efficient tool.
In the context of solvable lattice models \cite{B}, 
it produces eigenvectors of transfer matrices from solutions of 
Bethe equations.
Beside exact evaluation of physical quantities,
Bethe ansatz has brought a number of applications also 
in representation theory and combinatorics.
A prominent example is the fermionic formula \cite{KKR,KR}, which 
grew out of the completeness problem, a certain 
counting of Bethe vectors under string hypotheses.
With the advent of the corner transfer matrix \cite{B} and 
the crystal base theory \cite{Ka},  the fermionic formulas 
have been reformulated as the so called  
$X=M$ conjecture for any affine Lie algebra 
\cite{HKOTY1, HKOTT}.
See \cite{Schi} for the current status of the conjecture.

The proof of the fermionic formula for $A^{(1)}_n$ \cite{KKR, KR} 
may be viewed as a combinatorial version of Bethe ansatz.
As a substitute of solutions to Bethe equations,
the combinatorial object called 
{\em rigged configuration} (RC) is introduced.
It is an $n$-tuple of Young diagrams (configuration) 
each row of which is assigned with an integer (rigging) obeying 
a selection rule.
The Bethe vectors are replaced by 
Littlewood-Richardson tableaux, or equivalently, 
{\em highest paths}.
The latter are $A_n$ highest weight elements in 
$B_{\mu_1} \otimes \cdots \otimes B_{\mu_m}$, where 
$B_l$ is the $A^{(1)}_n$ crystal of the $l$-fold symmetric tensor 
representation corresponding to a ``local spin". 
Under these setting, Bethe ansatz should produce  
highest paths from rigged configurations and vice versa.
What achieves this and thereby proves the fermionic formula 
is the celebrated Kerov-Kirillov-Reshetikhin (KKR) bijection.
It is defined by a purely combinatorial algorithm 
on rigged configurations and highest paths, whose meaning 
however has remained rather mysterious.

The purpose of this paper is to clarify 
the representation theoretical origin of the KKR bijection
in the light of crystal base theory and 
associated soliton cellular automata \cite{HKOTY2}.
Let $p$ be the highest path that corresponds to a rigged configuration 
under the KKR bijection.
Our main Theorem \ref{th:main} asserts that 
\begin{equation}\label{eq:p0}
p = \Phi_1{\mathcal C}_1\Phi_{2}{\mathcal C}_{2}
\cdots \Phi_n{\mathcal C}_n(p^{(n)}).
\end{equation}
Here $p^{(n)}$ is a trivial vacuum path and 
$\Phi_1{\mathcal C}_1\Phi_{2}{\mathcal C}_{2}
\cdots \Phi_n{\mathcal C}_n$ is a vertex operator 
in the sense explained below.
Recall that the $A^{(1)}_n$ crystal 
$B_l$ consists of length $l$ row semistandard tableaux 
with letters $1,2, \ldots, n+1$.
As a set the affine crystal is given by 
$\Aff(B_l) = \{b[d] \mid b \in B_l, d \in \Z \}$ 
by assigning the mode $d$ to $B_l$.
The isomorphism 
$\Aff(B_l) \otimes \Aff(B_m) 
\simeq \Aff(B_m) \otimes \Aff(B_l)$
is called the {\em combinatorial $R$} \cite{KMN,NY}.
We will deal with 
the nested family of algebras 
$A^{(1)}_n \supset A^{(1)}_{n-1} \supset \cdots \supset A^{(1)}_0$ 
and the associated crystals:
\begin{equation}\label{eq:nest}
B_l = B^{\ge 1}_l \supset B^{\ge 2}_l \supset
\cdots \supset B^{\ge n+1}_l,
\end{equation}
where the superscript means the restriction on tableau letters. 
Given a rigged configuration (\ref{eq:rc}),
the right hand side of (\ref{eq:p0}) is determined 
by using the combinatorial $R$ for the crystals (\ref{eq:nest}) alone.
In fact $p^{(a)} = \Phi_{a+1}{\mathcal C}_{a+1}
\cdots \Phi_n{\mathcal C}_n(p^{(n)})$
becomes a highest path 
in $B^{\ge a+1}_{\mu_1} \otimes \cdots \otimes 
B^{\ge a+1}_{\mu_m}$
with respect to $A_{n-a}$.
Here $\mu=\mu^{(a)}$ is the $a$-th Young diagram 
in the rigged configuration. 
In particular $p^{(n)}$ is trivially fixed.
The map ${\mathcal C}_a$ sends $p^{(a)}$ to $\Aff(B^{\ge a+1}_{\mu_1}) 
\otimes \cdots \otimes \Aff(B^{\ge a+1}_{\mu_m})$
by assigning modes based on the rigging 
with certain normal ordering afterwards.
Then the map $\Phi_a$ creates a highest path $p^{(a-1)}$ in 
$B^{\ge a}_{\lambda_1} \otimes \cdots \otimes B^{\ge a}_{\lambda_l}$
$(\lambda = \mu^{(a-1)})$
by using the data ${\mathcal C}_a(p^{(a)})$ and 
the natural embedding 
$B^{\ge a+1}_l \hookrightarrow B^{\ge a}_l$ as sets.
Thus the composition 
$\Phi_1{\mathcal C}_1\Phi_{2}{\mathcal C}_{2}
\cdots \Phi_n{\mathcal C}_n$ grows 
the trivial $A_0$ path $p^{(n)}$ into an 
$A_n$ highest path by gradually taking the 
rigged configuration into account.
Such a usage of the family 
$A^{(1)}_n \supset A^{(1)}_{n-1} \supset \cdots \supset A^{(1)}_0$
is the typical strategy in the nested Bethe ansatz.
In fact our formula (\ref{eq:p0}) is 
a crystal theoretical formulation, i.e., 
$q=0$ analogue, of Schultz's 
construction of Bethe vectors \cite{Schu}.

The result (\ref{eq:p0}) admits a further interpretation 
in the realm of soliton theory, which we shall now explain.
Consider the crystal 
$B_{\mu_1} \otimes \cdots \otimes B_{\mu_m}$, 
where we formally take $m \rightarrow \infty$ 
and impose the boundary condition 
$p_j = \boxed{1\ldots 1} \in B_{\mu_j}$ 
for $j \gg 1$ on its elements $p=p_1 \otimes p_2 \otimes \cdots$. 
Then the system can be 
endowed with a commuting family of 
time evolutions and behaves as 
a soliton cellular automaton.
It was invented as the box-ball system \cite{TS,T}
and subsequently reformulated 
by the crystal theory \cite{HHIKTT,FOY,HKOTY2}.
Here is a typical time evolution pattern when $\mu_j = 1$ for all $j$.
(We omit $\otimes$ and write $\boxed{1}$ simply as $1$, etc.)
 
\noindent
\begin{center}
$ t=0: \quad 1111222211113321111411111111111111111111111 \cdots $\\ 
$ t=1: \quad 1111111122221113321141111111111111111111111 \cdots $\\ 
$ t=2: \quad 1111111111112222113324111111111111111111111 \cdots $\\
$ t=3: \quad 1111111111111111222213432111111111111111111 \cdots $\\
$ t=4: \quad 1111111111111111111122321432211111111111111 \cdots $\\
$ t=5: \quad 1111111111111111111111213221143221111111111 \cdots $\\
$ t=6: \quad 1111111111111111111111121113221114322111111 \cdots $\\
$ t=7: \quad 1111111111111111111111112111113221111432211 \cdots $
\end{center}

\noindent
The three solitons with amplitudes 4, 3 and 1 regain
the original amplitudes after the collision while interchanging 
the internal degrees of freedom.
The dynamics is governed by the combinatorial $R$ and 
as a result, highest paths remain highest.
In the above example, all the paths are highest indeed.
(The paths regarded as words are lattice permutations.)
Thus it is natural to ask; 
what kind of time evolutions does the soliton cellular automaton 
induce on the rigged configurations?
Our answer to the question is summarized in the following table.

\begin{table}[h]
\begin{center}
\begin{tabular}{c|c|c}
\hfil Bethe ansatz & Crystal base theory&  Soliton cellular automaton \\
\hline 
&&\vspace{-0.2cm}\\
rigged configuration & 
$\Aff(B_{\lambda_1})\ot \cdots \ot \Aff(B_{\lambda_l})$ &
action-angle variable \\
&& \vspace{-0.3cm} \\
\hline
&& \vspace{-0.2cm}\\
KKR bijection & vertex operator 
& (inverse) scattering map
\end{tabular}
\end{center}
\end{table}

We show that under time evolutions 
the configuration remains unchanged while 
the rigging flows linearly.
In this sense the rigged configurations are action-angle variables
of the associated soliton cellular automaton.
The first Young diagram $\mu^{(1)}$ in the configuration 
gives the list of amplitudes of solitons.
For example, $t=4$ state in the above pattern coincides with  
$p$ in Example \ref{ex:okado} (apart from redundant $1$'s) 
for which $\mu^{(1)}=(4,3,1)$ indeed.
The meaning of the other Young diagrams 
is understood similarly along the nested family of 
cellular automata explained in Section \ref{subsec:spectroscopy}.

Finally to interpret the formula (\ref{eq:p0}),
write it as $p=\Phi_1{\mathcal C}_1(p^{(1)})$.
In the above example, one has 
$p^{(1)}= \boxed{2222}\otimes \boxed{233} \otimes \boxed{4}$,
and this is nothing but the list of incoming solitons.
The map ${\mathcal C}_1$ assigns $p^{(1)}$ with the modes 
in affine crystals that encode the positions of solitons.
Thus ${\mathcal C}_1(p^{(1)})$  
is the {\em scattering data} of the soliton cellular automaton.
Then $\Phi_1$ plays the role of a vertex operator to create 
solitons over the ``vacuum path" $111\ldots $ by injecting 
the scattering data ${\mathcal C}_1(p^{(1)})$ 
by using combinatorial $R$.
Again these constructions work
inductively along the nested family of crystals for 
$A^{(1)}_n \supset A^{(1)}_{n-1} \supset \cdots \supset A^{(1)}_0$.
Thus the composition 
$\Phi_1{\mathcal C}_1 \cdots \Phi_n{\mathcal C}_n$ 
is the inverse scattering map that reproduces solitons from
the scattering data. 

The layout of the paper is as follows.
In Section \ref{sec:A(1)}, we treat the $A^{(1)}_n$ case.
A proof of Theorem \ref{th:main} will be presented in \cite{S}.
Section \ref{subsec:spectroscopy} includes the 
solution of the direct scattering problem, namely, 
a method to determine the rigged configuration 
from a given highest path by only using combinatorial $R$ and 
the KKR bijection for $\widehat{sl}_2$. 
It is a crystal theoretical separation of variables.
In Section \ref{sec:conj}, 
conjectures similar to (\ref{eq:p0}) are presented 
for the bijection in the 
other non-exceptional affine algebras \cite{OSS,SS}.
They all possess a similar feature to the $A^{(1)}_n$ case, 
but the items as $\Phi_a, {\mathcal C}_a$ 
need individual descriptions.
We do this from Sections \ref{sec:D(1)} to \ref{sec:III}.
Appendix \ref{app:crystal} is a brief exposition 
of the crystal base theory.
Further applications of the present results to 
fermionic formulas and soliton cellular automata
will be given elsewhere.

\section{$A^{(1)}_n$ case}\label{sec:A(1)}

\subsection{Rigged configurations}\label{subsec:rcA(1)}
Consider the data of the form
\begin{equation}\label{eq:rc}
(\mu^{(0)}, (\mu^{(1)},J^{(1)}), \ldots, (\mu^{(n)},J^{(n)})),
\end{equation}
where 
$\mu^{(a)}=(\mu^{(a)}_1,\ldots, \mu^{(a)}_{l_a})$ is a partition 
and 
$J^{(a)}= (J^{(a)}_1, \ldots, J^{(a)}_{l_a}) \in (\Z_{\ge 0})^{l_a}$.
Set 
\begin{align}
E^{(a)}_j = \sum_{i=1}^{l_a}\min(j,\mu^{(a)}_i). \label{eq:Eaj}
\end{align}
Define the vacancy numbers by
\begin{equation}\label{eq:paj}
p^{(a)}_j = E^{(a-1)}_j-2E^{(a)}_j + E^{(a+1)}_j\quad (1 \le a \le n),
\end{equation}
where $E^{(n+1)}_j = 0$.
The data (\ref{eq:rc}) is called a rigged configuration 
if the following conditions are satisfied for all
$1 \le a \le n$ and $j \in \Z_{>0}$:
\begin{equation}\label{eq:cond}
0\le J^{(a)}_i \le J^{(a)}_{i+1}\le \cdots \le J^{(a)}_l 
\le p^{(a)}_j \; \hbox{ if  } \{i,i+1,\ldots, l\}
=\{k\mid \mu^{(a)}_k = j\}.
\end{equation}
The array of partitions 
$\mu^{(0)}, \ldots, \mu^{(n)}$ is called a configuration 
and the nonnegative integers $J^{(a)}_i$ are called rigging.
If $\mu^{(a)}$ is represented as a Young diagram,
the vacancy number $p^{(a)}_j$ is assigned to each 
``cliff" of width $j$.
$J^{(a)}$ may be viewed as assigning to the cliff 
a partition whose parts are at most $p^{(a)}_j$.
Here is an example of the 
$A^{(1)}_3$ rigged configuration.
\begin{equation}\label{eq:rcA(1)}
\unitlength 0.1in
\begin{picture}(  0.0000,  11.000)(  16.0000,-14.5000)
\put(7,-5){$\mu^{(0)}$}
\put(12.5,-5){$\mu^{(1)}$}
\put(17,-5){$\mu^{(2)}$}
\put(22,-5){$\mu^{(3)}$}
%
\special{pn 8}%
\special{pa 600 600}%
\special{pa 600 1320}%
\special{fp}%
%
\special{pn 8}%
\special{pa 720 600}%
\special{pa 720 1320}%
\special{fp}%
%
\special{pn 8}%
\special{pa 600 720}%
\special{pa 960 720}%
\special{fp}%
%
\special{pn 8}%
\special{pa 600 600}%
\special{pa 960 600}%
\special{fp}%
%
\special{pn 8}%
\special{pa 840 600}%
\special{pa 840 840}%
\special{fp}%
%
\special{pn 8}%
\special{pa 840 600}%
\special{pa 960 600}%
\special{pa 960 720}%
\special{pa 840 720}%
\special{pa 840 600}%
\special{fp}%
%
\special{pn 8}%
\special{pa 720 720}%
\special{pa 840 720}%
\special{pa 840 840}%
\special{pa 720 840}%
\special{pa 720 720}%
\special{fp}%
%
\special{pn 8}%
\special{pa 600 720}%
\special{pa 720 720}%
\special{pa 720 840}%
\special{pa 600 840}%
\special{pa 600 720}%
\special{fp}%
%
\special{pn 8}%
\special{pa 600 840}%
\special{pa 720 840}%
\special{pa 720 960}%
\special{pa 600 960}%
\special{pa 600 840}%
\special{fp}%
%
\special{pn 8}%
\special{pa 600 960}%
\special{pa 720 960}%
\special{pa 720 1080}%
\special{pa 600 1080}%
\special{pa 600 960}%
\special{fp}%
%
\special{pn 8}%
\special{pa 600 1080}%
\special{pa 720 1080}%
\special{pa 720 1200}%
\special{pa 600 1200}%
\special{pa 600 1080}%
\special{fp}%
%
\special{pn 8}%
\special{pa 600 1200}%
\special{pa 720 1200}%
\special{pa 720 1320}%
\special{pa 600 1320}%
\special{pa 600 1200}%
\special{fp}%
%
\special{pn 8}%
\special{pa 600 1320}%
\special{pa 720 1320}%
\special{pa 720 1440}%
\special{pa 600 1440}%
\special{pa 600 1320}%
\special{fp}%
%
\special{pn 8}%
\special{pa 1440 600}%
\special{pa 1440 840}%
\special{fp}%
%
\special{pn 8}%
\special{pa 1200 600}%
\special{pa 1200 1080}%
\special{fp}%
%
\special{pn 8}%
\special{pa 1320 600}%
\special{pa 1320 1080}%
\special{fp}%
%
\special{pn 8}%
\special{pa 1200 600}%
\special{pa 1440 600}%
\special{fp}%
%
\special{pn 8}%
\special{pa 1200 720}%
\special{pa 1440 720}%
\special{fp}%
%
\special{pn 8}%
\special{pa 1200 840}%
\special{pa 1440 840}%
\special{fp}%
%
\special{pn 8}%
\special{pa 1200 960}%
\special{pa 1320 960}%
\special{fp}%
%
\special{pn 8}%
\special{pa 1200 1080}%
\special{pa 1320 1080}%
\special{fp}%
%
\special{pn 8}%
\special{pa 1680 600}%
\special{pa 1680 840}%
\special{fp}%
%
\special{pn 8}%
\special{pa 1800 600}%
\special{pa 1800 840}%
\special{fp}%
%
\special{pn 8}%
\special{pa 1680 600}%
\special{pa 1920 600}%
\special{fp}%
%
\special{pn 8}%
\special{pa 1680 720}%
\special{pa 1920 720}%
\special{fp}%
%
\special{pn 8}%
\special{pa 1920 600}%
\special{pa 1920 720}%
\special{fp}%
%
\special{pn 8}%
\special{pa 1680 840}%
\special{pa 1800 840}%
\special{fp}%
%
\special{pn 8}%
\special{pa 2160 600}%
\special{pa 2280 600}%
\special{pa 2280 720}%
\special{pa 2160 720}%
\special{pa 2160 600}%
\special{fp}%
\put(11.8200,-6.9000){\makebox(0,0)[rt]{0}}%
\put(11.7600,-9.3000){\makebox(0,0)[rt]{1}}%
\put(16.6800,-6.1800){\makebox(0,0)[rt]{1}}%
\put(16.7400,-7.6200){\makebox(0,0)[rt]{1}}%
\put(21.4200,-6.1800){\makebox(0,0)[rt]{0}}%
\put(18.1800,-7.6200){\makebox(0,0)[lt]{0}}%
\put(13.3800,-9.9000){\makebox(0,0)[lt]{1}}%
\put(19.4400,-6.3000){\makebox(0,0)[lt]{1}}%
\put(23.0400,-6.2400){\makebox(0,0)[lt]{0}}%
\put(13.5000,-8.7600){\makebox(0,0)[lt]{0}}%
\put(14.5800,-6.3000){\makebox(0,0)[lt]{0}}%
\put(14.5800,-7.5600){\makebox(0,0)[lt]{0}}%
\end{picture}%
\end{equation}

\noindent
The vacancy number is written on the left of the 
Young diagrams and the rigging is attached to each row. 
For a partition $\lambda=(\lambda_1,\ldots, \lambda_k)$, 
we let ${\rm RC}(\lambda)$ 
denote the set of rigged configurations 
(\ref{eq:rc}) with $\mu^{(0)}=\lambda$.

\subsection{\mathversion{bold} Crystals}
We recapitulate basic facts on 
the $A^{(1)}_n$ crystal $B_l$.
For a general background see Appendix \ref{app:crystal}.
The $B_l$ is the crystal base of the $l$-fold symmetric 
tensor representation. 
As the set it is given by 
\begin{equation}\label{eq:Bl}
B_l = \{x=(x_1,\ldots, x_{n+1})\in (\Z_{\ge 0})^{n+1}
\mid x_1+ \cdots + x_{n+1}=l\}.
\end{equation}
The Kashiwara operators act as 
$\tilde{e}_i(x)=x', \tilde{f}_i(x)=x''$ with
$x'_j = x_j + \delta_{i,j}-\delta_{i,j+1}$ and 
$x''_j = x_j - \delta_{i,j}+\delta_{i,j+1}$.
Here indices are in $\Z_{n+1}$ 
and $x'$ and $x''$ are to be understood as 0 
unless they belong to $(\Z_{\ge 0})^{n+1}$.
The combinatorial $R: 
\Aff(B_l) \ot \Aff(B_m) \rightarrow 
\Aff(B_m) \ot \Aff(B_l)$ has the form 
$R: x[d]\otimes y[e] \mapsto
\tilde{y}[e-H(x\ot y)]\ot \tilde{x}[d+H(x\ot y)]$ with 
\begin{align}
&{\tilde x}_i = x_i+Q_i(x,y)-Q_{i-1}(x,y),\quad 
{\tilde y}_i = y_i+Q_{i-1}(x,y)-Q_i(x,y),\\
&Q_i(x,y) = \min \{ \sum_{j=1}^{k-1}x_{i+j} + \sum_{j=k+1}^{n+1} y_{i+j} 
\mid 1 \le k \le n+1 \}, \\
&H(x\ot y) = \min(l,m)-Q_0(x,y).\label{eq:h}
\end{align}
The energy function $H$ here is 
normalized so that 
$0 \le H \le \min(l,m)$ and 
coincides with the ``winding number" \cite{NY}.
The element $x=(x_1,\ldots, x_{n+1})$ is also
denoted by a row shape semistandard tableau
of length $l$ containing the letter $i$ $x_i$ times and 
$x[d]\in \Aff(B_l)$ by
the tableau with index $d$.
For example in $A^{(1)}_3$, the following 
stand for the same relation under $R$:
\begin{align*}
(1,2,0,1)[5]\ot (1,0,1,0)[9] &\simeq (0,1,0,1)[8]\ot (2,1,1,0)[6],\\
\fbox{1224}_{\,5}\ot \fbox{13}_{\,9} &\simeq 
\fbox{24}_{\,8}\ot \fbox{1123}_{\,6}.
\end{align*}
To save the space we use the notation:
\begin{equation}\label{eq:al}
a^l = \boxed{a\cdots a} \in B_l.
\end{equation}
The relation $x\otimes y \simeq {\tilde y} \otimes {\tilde x}$
is depicted as

\begin{equation}\label{eq:Rfig}
\begin{picture}(20,18)(20,7)
\put(0,10){\line(1,0){20}}\put(-6,8){$x$}\put(22,8){${\tilde x}$}
\put(10,0){\line(0,1){20}}\put(8,25){$y$}\put(8,-9){${\tilde y}$}
\put(47,8){or}
\put(80,10){\line(1,0){20}}\put(74,8){${\tilde y}$}\put(102,8){$y$}
\put(90,0){\line(0,1){20}}\put(88,25){${\tilde x}$}\put(88,-9){$x$}
\end{picture}
\end{equation}

\vspace{0.7cm}
Setting 
\begin{equation}
B^{\ge a+1}_l = 
\{(x_1,\ldots, x_{n+1}) \in B_l \mid x_1=\cdots =x_a=0\}\quad 
(0 \le a \le n),
\end{equation}
we have 
\begin{equation}\label{eq:embed}
B_l = B^{\ge 1}_l \supset B^{\ge 2}_l \supset
\cdots \supset B^{\ge n+1}_l =\{(n+1)^l\} 
\end{equation}
as sets.
We will need to consider the crystals not only 
for $A^{(1)}_n$ but also for the nested family
$A^{(1)}_0, A^{(1)}_1, \ldots, A^{(1)}_{n-1}$.
In such a circumstance we realize 
the crystal $B_l$ for $A^{(1)}_{n-a}\; (0 \le a \le n)$ 
on the set $B^{\ge a+1}_l$ with the Kashiwara operators 
$\tilde{e}_i, \tilde{f}_i\; (a \le i \le n)$.
In this convention the 
highest element with respect to $A_{n-a}$ is 
$(a+1)^l \in B^{\ge a+1}_l$.

Let 
\begin{equation}\label{eq:P}
{\mathcal P}_+(\lambda) = 
\{p \in B_{\lambda_1}\otimes \cdots \otimes B_{\lambda_k}
\mid {\tilde e}_ip=0, \; 1 \le i \le n\}.
\end{equation}
be the set of highest elements (paths) with respect to $A_n$.
The bijection \cite{KKR,KR} between 
${\rm RC}(\lambda)$ and the Littlewood-Richardson 
tableaux is translated to the one between ${\rm RC}(\lambda)$ and 
${\mathcal P}_+(\lambda)$.
We call the resulting map the KKR bijection.
See \cite{Schi} for a recent review.
It sends the rigged configuration (\ref{eq:rcA(1)}) to
\begin{equation}\label{eq:pA(1)}
\boxed{111}\otimes \boxed{22} \otimes \boxed{3}
\otimes \boxed{1} \otimes \boxed{4} \otimes \boxed{2}
\otimes \boxed{3}.
\end{equation}

\subsection{Normal ordering}\label{subsec:noA(1)}
For an element $b_1[d_1]\ot \cdots \ot b_m[d_m] 
\in \Aff(B_{l_1})\ot \cdots \ot \Aff(B_{l_m})$ 
we call the number $d_i$ the $i$-th mode.
By using the combinatorial $R$, 
tensor products can be reordered and the 
modes are changed accordingly.
Given an element $s \in \Aff(B_{l_1})\ot \cdots \ot \Aff(B_{l_m})$, 
define  ${\mathcal S}_m$ to be the set of such reordering as
\begin{equation}\label{eq:sm}
{\mathcal S}_m = \{s' \in \bigsqcup_{\sigma \in {\mathfrak S}_m}
\!\!'\;\Aff(B_{l_{\sigma(1)}})\ot \cdots \ot \Aff(B_{l_{\sigma(m)}}) \mid 
s' \simeq s \},
\end{equation}
where $\sqcup'$ means the disjoint union over $\sigma \in {\mathfrak S}_m$
satisfying $\sigma(i)<\sigma(j)$ for any $i,j$ such that $i<j$ and
$l_i=l_j$.
The cardinality of ${\mathcal S}_m$ is $m!$ if 
$l_1,\ldots, l_m$ are distinct.
For $i=2,\ldots, m$, let ${\mathcal S}_{i-1}$ 
be the subset of ${\mathcal S}_{i}$ having the maximal
$i$-th mode.
Then we have 
\begin{equation}
\emptyset \neq {\mathcal S}_1 \subseteq {\mathcal S}_2 \subseteq 
\cdots \subseteq {\mathcal S}_m.
\end{equation}
We call the elements of ${\mathcal S}_1$ 
{\em normal ordered forms} of $s$.
In general the normal ordered form 
$b_1[d_1] \ot \cdots \ot b_m[d_m]$ is not unique but 
the mode sequence $d_1,\ldots, d_m$ is unique by the definition.
For type $A^{(1)}_n$, one has $d_1 \le \cdots \le d_m$.
Any element of ${\mathcal S}_1$ is denoted by $:\!s\!:$.

In the context of soliton cellular automaton explained 
in Section \ref{sec:1}, 
the restriction on the union $\sqcup'$ in (\ref{eq:sm}) 
reflects the fact that solitons of equal velocity(=amplitude)
do not collide with each other.

\begin{example}
Take $s=\fbox{$1223$}_{\,5}\ot \fbox{$34$}_{\,8} 
\ot \fbox{$1$}_{\,6}$ for $A^{(1)}_{3}$.
\begin{align*}
{\mathcal S}_3 &= \{
\fbox{$1223$}_{\,5}\ot \fbox{$34$}_{\,8} \ot \fbox{$1$}_{\,6},\;\;
\fbox{$23$}_{\,8}\ot \fbox{$1234$}_{\,5} \ot \fbox{$1$}_{\,6},\;\;
\fbox{$23$}_{\,8}\ot \fbox{$4$}_{\,5} \ot \fbox{$1123$}_{\,6},\\
&\qquad \fbox{$3$}_{\,5}\ot \fbox{$24$}_{\,8} \ot 
\fbox{$1123$}_{\,6},\;\;
\fbox{$3$}_{\,5}\ot \fbox{$1224$}_{\,5} \ot \fbox{$13$}_{\,9},\;\;
\fbox{$1223$}_{\,5}\ot \fbox{$4$}_{\,5} \ot \fbox{$13$}_{\,9}\},\\
{\mathcal S}_2&={\mathcal S}_1 = \{
\fbox{$3$}_{\,5}\ot \fbox{$1224$}_{\,5} \ot \fbox{$13$}_{\,9},\;\;
\fbox{$1223$}_{\,5}\ot \fbox{$4$}_{\,5} \ot \fbox{$13$}_{\,9}\}.
\end{align*}
\end{example}

\subsection{\mathversion{bold} Maps ${\mathcal C}_1, \ldots, {\mathcal C}_n$}
\label{subsec:clclA(1)}
Pick the color $a$ part $(\mu^{(a)},J^{(a)})$ of the rigged 
configuration.  Here we simply write it as $(\mu, J)$. 
Namely $\mu=(\mu_1, \ldots, \mu_m)$ is a partition and 
$J=(J_i)$, where $J_i$ is the rigging attached to the $i$-th row 
in $\mu$ of length $\mu_i$.
For $1 \le a \le n$, let $B_l = B^{\ge a+1}_l$ be the 
$A^{(1)}_{n-a}$ crystal in the sense explained around 
(\ref{eq:embed}). 
Define the map ${\mathcal C}_a$ among the $A^{(1)}_{n-a}$ crystals by
\begin{align}
{\mathcal C}_a: &\;B_{\mu_1} \ot \cdots \ot B_{\mu_m} \rightarrow \;
:\Aff(B_{\mu_1}) \ot \cdots \ot \Aff(B_{\mu_m}):
\quad (1\le a \le n)\nonumber \\
&\quad b_1\ot \cdots \ot b_m \quad \mapsto\quad 
:\!b_1[d_1]\ot \cdots \ot b_m[d_m]\!: \label{eq:Ca}\\
d_i &= J_i + \sum_{0 \le k < i}H(b_k\ot b^{(k+1)}_i), \quad 
b_0=(a+1)^{\mu_1}.\label{eq:di}
\end{align}
Here $b^{(j)}_i \in B_{\mu_i}\; (j \le i)$ is defined by 
bringing $b_i$ to the left by the combinatorial $R$ as
\begin{equation}
(b_{j}\ot \cdots \ot b_{i-1}) \ot b_i \simeq 
b^{(j)}_i \ot ( \,\cdots)
\end{equation}
under the isomorphism 
$(B_{\mu_j}\ot \cdots \ot B_{\mu_{i-1}})\ot B_{\mu_i}
\simeq B_{\mu_i} \ot (B_{\mu_j}\ot \cdots \ot B_{\mu_{i-1}})$.

The map ${\mathcal C}_n$ involves  ``$A^{(1)}_0$ crystal" 
$B^{\ge n+1}_l=\{ (n+1)^l\}$.
The following suffices to define ${\mathcal C}_n$:
\begin{equation}\label{eq:A0}
(n+1)^l \otimes (n+1)^m \simeq (n+1)^m \otimes (n+1)^l,
\quad 
H((n+1)^l \otimes (n+1)^m) = \min(l,m).
\end{equation}

Since the normal ordering in (\ref{eq:Ca}) is not unique,
${\mathcal C}_a$ is actually multi-valued in general.
Here we mean by ${\mathcal C}_a(\cdot)$ 
to pick any one of the normal ordered forms.

\subsection{\mathversion{bold} Maps $\Phi_1, \ldots, \Phi_n$}
\label{subsec:Phi}
Pick the color $a$ and $a-1$ parts of the configuration
and denote them simply by 
$\mu^{(a)}=(\mu_1,\ldots, \mu_m)$ and 
$\mu^{(a-1)}=(\lambda_1,\ldots, \lambda_k)$.
Set $B_l = B^{\ge a+1}_l$ and 
$B'_l = B^{\ge a}_l$.
We define the map $\Phi_a$ from the 
normal ordered elements in 
$A^{(1)}_{n-a}$ affine crystals to 
$A^{(1)}_{n-a+1}$ crystals:
\begin{align}
\Phi_a: \; :\Aff(B_{\mu_1})\ot \cdots \ot 
\Aff(B_{\mu_m}):
&\rightarrow 
B'_{\lambda_1} \ot \cdots \ot B'_{\lambda_k}
\quad (1 \le a \le n)\nonumber \\
b_1[d_1]\ot \cdots \ot b_m[d_m] \quad \quad\;\;
&\mapsto\quad 
c_1 \ot \cdots \ot c_k.\label{eq:Phi}
\end{align}
Since $b_1[d_1]\ot \cdots \ot b_m[d_m]$ is normal ordered, 
we know that $d_1 \le \cdots \le d_m$. 
We suppose further that $d_1 \ge 0$, which is the case 
in the actual use later.
Then the image $c_1 \ot \cdots \ot c_k$ is determined by
the following relation under the isomorphism of
$A^{(1)}_{n-a+1}$ crystals:
(We write ${\mathcal T}_a^{d} 
= \boxed{a}^{\ot d} \in (B^{\ge a}_1)^{\ot d}$ 
for short.)
\begin{equation}
\begin{split}
&({\mathcal T}_a^{d_1}\ot b_1 \ot {\mathcal T}_a^{d_2-d_1} \ot b_2 \ot \cdots 
\ot {\mathcal T}_a^{d_m-d_{m-1}}\ot b_m)\ot 
(a^{\lambda_1} \ot a^{\lambda_2} \ot \cdots \ot a^{\lambda_k})\\
& \simeq 
(c_1\ot \cdots \ot c_k) \ot \hbox{tail},
\end{split}\label{eq:zdon}
\end{equation}
Here we are regarding $b_i \in B_{\mu_i}=B^{\ge a+1}_{\mu_i}$ as 
an element of $B'_{\mu_i}=B^{\ge a}_{\mu_i}$ by the natural 
embedding (\ref{eq:embed}) as sets.
The tail part has the same structure as 
$({\mathcal T}_a^{d_1}\ot b_1 \ot {\mathcal T}_a^{d_2-d_1} \ot \cdots \ot b_m)$ 
on the left hand side.
In the actual use, it turns out to be 
$({\mathcal T}_a^{d_1}\ot a^{\mu_1}\ot {\mathcal T}_a^{d_2-d_1} \ot  \cdots 
\ot a^{\mu_m})$ containing the letter $a$ only.
(This fact will not be used.)

To obtain $c_1\ot \cdots \ot c_k$ using (\ref{eq:zdon}),
one applies the combinatorial $R$ on  
$B'_j \ot B'_l$
to carry $({\mathcal T}_a^{d_1}\ot b_1 \ot  \cdots 
\ot {\mathcal T}_a^{d_m-d_{m-1}} \ot b_m)$ through 
$(a^{\lambda_1} \ot a^{\lambda_2} \ot \cdots \ot a^{\lambda_k})$
to the right.  
The procedure is depicted as 

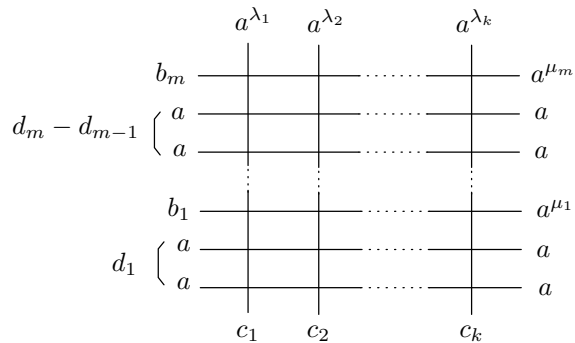
\begin{equation*}\unitlength 0.1in
\begin{picture}(  6.0000,  20.0000)( 47.8000, -26.2000)
\put(45.5000,-7.4000){\makebox(0,0){$a^{\lambda_1}$}}%
\put(60.3000,-19.3000){\makebox(0,0)[lt]{$a$}}%
\put(60.2000,-14.1000){\makebox(0,0)[lt]{$a$}}%
\put(60.2000,-12.1000){\makebox(0,0)[lt]{$a$}}%
\put(60.4000,-21.3000){\makebox(0,0)[lt]{$a$}}%
\put(60.0000,-9.8000){\makebox(0,0)[lt]{$a^{\mu_m}$}}%
\put(60.3000,-16.9000){\makebox(0,0)[lt]{$a^{\mu_1}$}}%
\put(56.9000,-23.9000){\makebox(0,0){$c_k$}}%
\put(48.9000,-23.9000){\makebox(0,0){$c_2$}}%
\put(45.2000,-23.9000){\makebox(0,0){$c_1$}}%
\put(39.4000,-19.7000){\makebox(0,0)[rt]{$d_1$}}%
\put(39.5000,-12.5000){\makebox(0,0)[rt]{$d_m-d_{m-1}$}}%
%
\special{pn 8}%
\special{pa 4520 1530}%
\special{pa 4520 1650}%
\special{dt 0.045}%
%
\special{pn 8}%
\special{pa 4520 1640}%
\special{pa 4520 2280}%
\special{fp}%
%
\special{pn 8}%
\special{pa 4520 870}%
\special{pa 4520 1510}%
\special{fp}%
%
\special{pn 8}%
\special{pa 5690 1540}%
\special{pa 5690 1660}%
\special{dt 0.045}%
%
\special{pn 8}%
\special{pa 5690 1650}%
\special{pa 5690 2290}%
\special{fp}%
%
\special{pn 8}%
\special{pa 5690 880}%
\special{pa 5690 1520}%
\special{fp}%
%
\special{pn 8}%
\special{pa 4890 1540}%
\special{pa 4890 1660}%
\special{dt 0.045}%
%
\special{pn 8}%
\special{pa 4890 1650}%
\special{pa 4890 2290}%
\special{fp}%
%
\special{pn 8}%
\special{pa 4890 880}%
\special{pa 4890 1520}%
\special{fp}%
\put(42.2000,-16.8000){\makebox(0,0)[rt]{$b_1$}}%
%
\special{pn 8}%
\special{pa 4260 1040}%
\special{pa 5100 1040}%
\special{fp}%
\special{pa 5460 1040}%
\special{pa 5940 1040}%
\special{fp}%
%
\special{pn 8}%
\special{pa 5100 1040}%
\special{pa 5460 1040}%
\special{dt 0.045}%
%
\special{pn 8}%
\special{pa 4260 1240}%
\special{pa 5100 1240}%
\special{fp}%
\special{pa 5460 1240}%
\special{pa 5940 1240}%
\special{fp}%
%
\special{pn 8}%
\special{pa 5100 1240}%
\special{pa 5460 1240}%
\special{dt 0.045}%
%
\special{pn 8}%
\special{pa 4260 1440}%
\special{pa 5100 1440}%
\special{fp}%
\special{pa 5460 1440}%
\special{pa 5940 1440}%
\special{fp}%
%
\special{pn 8}%
\special{pa 5100 1440}%
\special{pa 5460 1440}%
\special{dt 0.045}%
%
\special{pn 8}%
\special{pa 4270 1750}%
\special{pa 5110 1750}%
\special{fp}%
\special{pa 5470 1750}%
\special{pa 5950 1750}%
\special{fp}%
%
\special{pn 8}%
\special{pa 5110 1750}%
\special{pa 5470 1750}%
\special{dt 0.045}%
%
\special{pn 8}%
\special{pa 4270 1950}%
\special{pa 5110 1950}%
\special{fp}%
\special{pa 5470 1950}%
\special{pa 5950 1950}%
\special{fp}%
%
\special{pn 8}%
\special{pa 5110 1950}%
\special{pa 5470 1950}%
\special{dt 0.045}%
%
\special{pn 8}%
\special{pa 4270 2150}%
\special{pa 5110 2150}%
\special{fp}%
\special{pa 5470 2150}%
\special{pa 5950 2150}%
\special{fp}%
%
\special{pn 8}%
\special{pa 5110 2150}%
\special{pa 5470 2150}%
\special{dt 0.045}%
\put(42.0000,-14.0000){\makebox(0,0)[rt]{$a$}}%
\put(41.9000,-12.1000){\makebox(0,0)[rt]{$a$}}%
\put(42.2000,-19.0000){\makebox(0,0)[rt]{$a$}}%
\put(42.2000,-21.0000){\makebox(0,0)[rt]{$a$}}%
\put(42.0000,-9.7000){\makebox(0,0)[rt]{$b_m$}}%
\put(57.1000,-7.5000){\makebox(0,0){$a^{\lambda_k}$}}%
\put(49.3000,-7.5000){\makebox(0,0){$a^{\lambda_2}$}}%
%
\special{pn 8}%
\special{pa 4060 1230}%
\special{pa 4030 1268}%
\special{pa 4030 1418}%
\special{pa 4060 1456}%
\special{pa 4060 1456}%
\special{fp}%
%
\special{pn 8}%
\special{pa 4080 1910}%
\special{pa 4050 1948}%
\special{pa 4050 2098}%
\special{pa 4080 2136}%
\special{pa 4080 2136}%
\special{fp}%
\end{picture}%
\end{equation*}

\subsection{Main theorem}\label{subsec:main}

Define the $A^{(1)}_0$ crystal element 
\begin{equation}\label{eq:pn}
p^{(n)} = (n+1)^{\mu^{(n)}_1}\ot \cdots \ot (n+1)^{\mu^{(n)}_{l_n}}.
\end{equation}
\begin{theorem}\label{th:main}
The image $p$ of the rigged configuration 
$(\mu^{(0)}, (\mu^{(1)},J^{(1)}), \ldots, (\mu^{(n)},J^{(n)}))$ 
under the KKR bijection is given by
\begin{equation}\label{eq:p}
p = \Phi_1{\mathcal C}_1\Phi_{2}{\mathcal C}_{2}\cdots \Phi_n{\mathcal C}_n(p^{(n)}).
\end{equation}
\end{theorem}
The theorem asserts that 
$\Phi_1{\mathcal C}_1\Phi_{2}{\mathcal C}_{2}\cdots \Phi_n{\mathcal C}_n(p^{(n)})$ 
is independent of the choices of the possibly 
non-unique normal ordered forms
when applying the maps ${\mathcal C}_1, \ldots, {\mathcal C}_n$.
A proof is presented in \cite{S}.

Set
\begin{equation}\label{eq:pa}
p^{(a)} = \Phi_{a+1}{\mathcal C}_{a+1}\cdots \Phi_n{\mathcal C}_n(p^{(n)})\quad 
(0 \le a \le n-1),
\end{equation}
which belongs to the $A^{(1)}_{n-a}$ crystal
$B^{\ge a+1}_{\mu^{(a)}_1} \ot \cdots 
\ot B^{\ge a+1}_{\mu^{(a)}_{l_a}}$.
Thus $p$ in (\ref{eq:p}) is $p^{(0)}$.
\begin{corollary}
For $0 \le a \le n-1$, 
$p^{(a)}$ coincides with the image of the 
truncated rigged configuration 
$(\mu^{(a)}, (\mu^{(a+1)},J^{(a+1)}), \ldots, (\mu^{(n)},J^{(n)}))$ 
under the KKR bijection.
\end{corollary}

\begin{remark}\label{rem:comp}
In a rigged configuration (\ref{eq:rc}), 
one may treat $\mu^{(a)}$ as a composition rather than a partition.
In that case, the condition (\ref{eq:cond}) should read 
$0 \le J^{(a)}_{i_1} \le \cdots \le J^{(a)}_{i_m} \le p^{(a)}_j$ 
if $\{i_1 < \cdots < i_m\} = \{k\mid \mu^{(a)}_k=j\}$.
The other necessary change is only to redefine $b_0$ in (\ref{eq:di})
as $b_0 = (a+1)^{\max\{\mu_1,\ldots, \mu_m\}}$.
Actually, $b_0 = (a+1)^M$ with any 
$M \ge \max\{\mu_1,\ldots, \mu_m\}$ leads to the same $d_i$
in (\ref{eq:di}).
The same remark applies also for the other algebras 
treated in this paper.
\end{remark}

\begin{example}
Consider the $A^{(1)}_3$ rigged configuration (\ref{eq:rcA(1)}).
According to (\ref{eq:pn}) we set 
\begin{equation*}
p^{(3)} = \boxed{4},\quad {\mathcal C}_3(p^{(3)}) = 
\boxed{4}_{\, 1}.
\end{equation*}
We use (\ref{eq:zdon}) to find 
$p^{(2)}=\Phi_3{\mathcal C}_3(p^{(3)})$.
It amounts to calculating 
${\mathcal T}^1_3 \otimes \boxed{4} \otimes 
(\; \boxed{33} \otimes \boxed{3}\;)$.
This is done as
\begin{equation*}
\unitlength 0.8mm
{\small
\begin{picture}(80,29)(-25,-19)
\multiput(0,-12)(0,12){2}{
\multiput(0,0)(18,0){2}{\line(1,0){10}}
\multiput(5,-3)(18,0){2}{\line(0,1){6}}
}

\put(3,5){33}\put(22,5){3}

\put(-4,-1){4}\put(13.5,-1){3}
\put(29,-1){3}

\put(3,-7){34}\put(22,-7){3}

\put(-4,-13){3}\put(13.5,-13){4}\put(29,-13){3}

\put(3,-19){33}\put(22,-19){4}

\end{picture}
}
\end{equation*}
Thus we have 
\begin{equation*}
p^{(2)} = \boxed{33} \otimes \boxed{4}, \quad 
{\mathcal C}_2(p^{(2)}) = 
\boxed{3}_{\,1} \otimes \boxed{34}_{\,3},
\end{equation*}
where we have used 
$H(\boxed{33} \otimes \boxed{4})=0, 
H(\boxed{33} \otimes \boxed{33})=2$ and 
$H(\boxed{33} \otimes \boxed{3})=1$.  
To find $p^{(1)}=\Phi_1{\mathcal C}_1(p^{(2)})$, 
we calculate 
${\mathcal T}^1_2 \otimes \boxed{3} \otimes 
{\mathcal T}^2_2 \otimes \boxed{34} \otimes
(\,\boxed{22} \otimes \boxed{22} \otimes \boxed{2}
\otimes \boxed{2} \,)$.
\begin{equation*}
\unitlength 0.8mm
{\small
\begin{picture}(80,67)(-10,-22)
\multiput(0,-12)(0,12){5}{
\multiput(0,0)(18,0){4}{\line(1,0){10}}
\multiput(5,-3)(18,0){4}{\line(0,1){6}}
}

\put(3,41){22}\put(21,41){22}
\put(40,41){2}\put(58,41){2}

\put(-6,35){34}\put(12,35){22}\put(30,35){22}
\put(48,35){22}\put(66,35){22}

\put(3,29){34}\put(21,29){22}
\put(40,29){2}\put(58,29){2}

\put(-4.6,23){2}\put(13.4,23){3}\put(31.4,23){2}
\put(49.4,23){2}\put(67,23){2}

\put(3,17){24}\put(21,17){23}
\put(40,17){2}\put(58,17){2}

\put(-4.6,11){2}\put(13.4,11){4}\put(31.4,11){2}
\put(49.4,11){2}\put(67,11){2}

\put(3,5){22}\put(21,5){34}
\put(40,5){2}\put(58,5){2}

\put(-4.6,-1){3}\put(13.4,-1){2}\put(31.4,-1){3}
\put(49.4,-1){2}\put(67,-1){2}

\put(3,-7){23}\put(21,-7){24}
\put(40,-7){3}\put(58,-7){2}

\put(-4.6,-13){2}\put(13.4,-13){3}\put(31.4,-13){4}
\put(49.4,-13){3}\put(67,-13){2}

\put(3,-19){22}\put(21,-19){23}
\put(40,-19){4}\put(58,-19){3}

\end{picture}
}
\end{equation*}
Thus we find 
\begin{equation*}
p^{(1)} = \boxed{22}\otimes \boxed{23} 
\otimes \boxed{4} \otimes \boxed{3},\quad 
{\mathcal C}_1(p^{(1)}) = \, :\,
\boxed{22}_{\,2}\otimes \boxed{23} _{\,3}
\otimes \boxed{4}_{\,1} \otimes \boxed{3}_{\,3} :
\end{equation*}
There are three normal ordered forms 
$\boxed{2}_{\,1}\otimes \boxed{3}_{\,2} 
\otimes \boxed{22}_{\,3} \otimes \boxed{34}_{\,3}$,\;
$\boxed{2}_{\,1}\otimes \boxed{23}_{\,2} 
\otimes \boxed{2}_{\,3} \otimes \boxed{34}_{\,3}$
and
$\boxed{2}_{\,1}\otimes \boxed{23}_{\,2} 
\otimes \boxed{24}_{\,3} \otimes \boxed{3}_{\,3}$.
Any one of them can be chosen as ${\mathcal C}_1(p^{(1)})$.
We illustrate the derivation of 
$p=\Phi_1{\mathcal C}_1(p^{(1)})$ 
along the first one.
According to (\ref{eq:zdon}) we calculate
${\mathcal T}^1_1 \otimes \boxed{2} \otimes 
 {\mathcal T}^1_1 \otimes \boxed{3} \otimes 
 {\mathcal T}^1_1 \otimes 
\boxed{22} \otimes \boxed{34} \otimes 
(\boxed{111} \otimes \boxed{11} \otimes \boxed{1} 
\otimes \boxed{1}\otimes \boxed{1}\otimes \boxed{1}
\otimes \boxed{1})$.
\begin{equation*}
\unitlength 0.8mm
{\small
\begin{picture}(130,94)(-5,-22)
\multiput(0,-12)(0,12){7}{
\multiput(0,0)(18,0){7}{\line(1,0){10}}
\multiput(5,-3)(18,0){7}{\line(0,1){6}}
}

\put(2,65){111}\put(21,65){11}
\put(40,65){1}\put(58,65){1}
\put(76,65){1}\put(94,65){1}\put(112,65){1}

\put(-6,59){34}\put(12,59){11}\put(30,59){11}
\put(48,59){11}\put(66,59){11}\put(84,59){11}
\put(102,59){11}\put(120,59){11}

\put(2,53){134}\put(21,53){11}
\put(40,53){1}\put(58,53){1}
\put(76,53){1}\put(94,53){1}\put(112,53){1}

\put(-6,47){22}\put(12,47){34}\put(30,47){11}
\put(48,47){11}\put(66,47){11}\put(84,47){11}
\put(102,47){11}\put(120,47){11}

\put(2,41){122}\put(21,41){34}
\put(40,41){1}\put(58,41){1}
\put(76,41){1}\put(94,41){1}\put(112,41){1}

\put(-5,35){1}\put(13,35){2}\put(31,35){3}
\put(49,35){1}\put(67,35){1}\put(85,35){1}
\put(103,35){1}\put(121,35){1}

\put(2,29){112}\put(21,29){24}
\put(40,29){3}\put(58,29){1}
\put(76,29){1}\put(94,29){1}\put(112,29){1}

\put(-5,23){3}\put(13,23){1}\put(31,23){2}
\put(49,23){3}\put(67,23){1}\put(85,23){1}
\put(103,23){1}\put(121,23){1}

\put(2,17){123}\put(21,17){14}
\put(40,17){2}\put(58,17){3}
\put(76,17){1}\put(94,17){1}\put(112,17){1}

\put(-5,11){1}\put(13,11){2}\put(31,11){4}
\put(49,11){2}\put(67,11){3}\put(85,11){1}
\put(103,11){1}\put(121,11){1}

\put(2,5){113}\put(21,5){12}
\put(40,5){4}\put(58,5){2}
\put(76,5){3}\put(94,5){1}\put(112,5){1}

\put(-5,-1){2}\put(13,-1){3}\put(31,-1){1}
\put(49,-1){4}\put(67,-1){2}\put(85,-1){3}
\put(103,-1){1}\put(121,-1){1}

\put(2,-7){112}\put(21,-7){23}
\put(40,-7){1}\put(58,-7){4}
\put(76,-7){2}\put(94,-7){3}\put(112,-7){1}

\put(-5,-13){1}\put(13,-13){2}\put(31,-13){3}
\put(49,-13){1}\put(67,-13){4}\put(85,-13){2}
\put(103,-13){3}\put(121,-13){1}

\put(2,-19){111}\put(21,-19){22}
\put(40,-19){3}\put(58,-19){1}
\put(76,-19){4}\put(94,-19){2}\put(112,-19){3}

\end{picture}
}
\end{equation*}
The bottom line yields 
the path $p = \Phi_1{\mathcal C}_1(p^{(1)})$ 
in agreement with (\ref{eq:pA(1)}).
\end{example}

\subsection{Finding rigged configuration from highest path}
\label{subsec:spectroscopy}
In this subsection we consider a map from 
highest paths to the corresponding RCs by
only using
crystal isomorphisms and the KKR bijection for the $\widehat{sl}_2$ case.
The original idea for this has already arisen in \cite{Ta}.
To do this, we
need not only crystals $B_l$ but also $B^{2,l}$, two-row tableaux. It is
known that the set
of semistandard tableaux of $k\times l$ rectangular shape $B^{k,l}$ admits
the
$U'_q(A^{(1)}_n)$-crystal structure and there is a generalized KKR bijection
from
the highest elements in $B^{k_1,l_1}\ot\cdots\ot B^{k_m,l_m}$ to RCs
\cite{KSS}.
For a quick study of these matters we recommend \cite{Schi}. 
(Note that the order of the tensor products of crystals 
is opposite from ours.) 
In particular, the readers
are recommended to see Definition 4.5 to understand the operations
$\mbox{lh},\mbox{ls},\mbox{lb},
\delta,i,j$ used in the proofs below.

Let $u^{k,l}$ be the highest element of $B^{k,l}$. $u^{k,l}$ is the
$k\times l$ tableau
whose entries in the $j$-th row are all $j$ for $1\le j\le k$. For a
partition $\la=(\la_1,\ldots,
\la_k)$ we define 
$B_\la=B^{1,\la_1}\ot\cdots\ot B^{1,\la_k}$
and $B^{(2)}_\la=B^{2,\la_1}\ot\cdots\ot B^{2,\la_k}$.
In this subsection we do not always write $\mu^{(0)}$ for the RC even when
it corresponds to
a path in $\P_+(\mu^{(0)})$.

\begin{proposition} \label{prop:time evol}
If an element $p$ of $\P_+(\la)$ corresponds to the RC
$((\mu^{(1)},J^{(1)}),\ldots,
(\mu^{(n)},J^{(n)}))$ under the bijection, then $u^{k,l}\ot p$ corresponds
to
\[
((\mu^{(1)},J^{(1)}),\ldots,(\mu^{(k)},\hat{J}^{(k)}),\ldots,(\mu^{(n)},J^{(
n)})),
\]
where $\hat{J}^{(k)}_i=J^{(k)}_i+\min(l,\mu^{(k)}_i)$.
\end{proposition}

\begin{proof}
Let $p^{(a)}_j$ (resp. $\tilde{p}^{(a)}_j$) be the vacancy number for $p$
(resp. $u^{k,l}\ot p$).
Note that $\tilde{p}^{(a)}_j-p^{(a)}_j=\min(l,j)\delta_{a,k}$. Noticing this
fact the proof goes
by induction on $|\la|$ using operations $\mbox{ls},\mbox{lh},i,\delta$.
\end{proof}

\begin{theorem}\label{th:type2}
Let $p$ be an element of $\P_+(\la)$ and let
$((\mu^{(1)},J^{(1)}),\ldots(\mu^{(n)},J^{(n)}))$ be
the corresponding RC. We assume that $p$ has $(u^{1,1})^{\ot L}$ at the
right end with
sufficiently large $L$. Set $\mu=\mu^{(1)}=(\mu_1,\ldots,\mu_l)$.
Consider $(u^{2,\mu_1}\ot\cdots\ot u^{2,\mu_l})\ot p$ and switch the order
by the combinatorial R:
\[
\begin{array}{ccc}
B^{(2)}_\mu\ot B_\la&\longrightarrow&B_\la\ot B^{(2)}_\mu \\
(u^{2,\mu_1}\ot\cdots\ot u^{2,\mu_l})\ot
p&\longmapsto&\tilde{p}\ot(b_1\ot\cdots\ot b_l).
\end{array}
\]
Then we have the following.
\begin{itemize}
\item[(1)] $\tilde{p}$ is an element of $\P_+(\la)$ that does not contain
letters greater than 2.
\item[(2)] $\tilde{p}$ corresponds to the RC $(\mu^{(1)},J^{(1)})$.
\item[(3)] The letters in the first row of $b_j$ are all 1 for any $j$.
\item[(4)] $\bar{b}_1\ot\cdots \ot \bar{b}_l$ is an element of $\P_+(\mu)$ that
corresponds to the RC 
$((\mu^{(2)},J^{(2)}),\ldots, (\mu^{(n)},J^{(n)}))$. 
Here $\bar{b}_j$ is the second row of $b_j$.
\end{itemize}
\end{theorem}

\begin{proof}
Let $RC=((\mu^{(1)},J^{(1)}),(\mu^{(2)},J^{(2)}),
\ldots,(\mu^{(n)},J^{(n)}))$
be the RC
corresponding to $p$. From Proposition \ref{prop:time evol} the RC
corresponding to
$(u^{2,\mu_1}\ot\cdots\ot u^{2,\mu_l})\ot p$ is given by
\[
\widehat{RC}=((\mu^{(1)},J^{(1)}),(\mu^{(2)},\hat{J}^{(2)}),
\ldots, (\mu^{(n)},J^{(n)}))
\]
where $\hat{J}^{(2)}_i=J^{(2)}_i+\sum_{k=1}^l\min(\mu_k,\mu^{(2)}_i)$. Since
the RC does not
change by the application of the combinatorial $R$ 
(see \cite{KSS} Lemma 8.5),
the RC corresponding
to $\tilde{p}\ot(b_1\ot\cdots\ot b_l)$ is also given by $\widehat{RC}$. Note
the following
facts.
\begin{itemize}
\item[(i)] A string in $\widehat{RC}$ is singular, if and only if the
corresponding string
	in $RC$ is singular.
\item[(ii)] There is no singular string in $(\mu^{(1)},J^{(1)})$.
\end{itemize}
We now apply the procedures $\mbox{ls},\mbox{lb},\mbox{lh},i,j,\delta$
successively on RC
to obtain $b_1\ot\cdots\ot b_l$ and $\tilde{p}$. 
By applying $\mbox{ls},i$, 
one can assume
$b_l\in B^{2,1}$. By applying $\mbox{lb},j,\mbox{lh},\delta$, we then obtain
$b_l=\boxed{\scriptstyle{1}  \atop a}$ ($a\ge2$). Notice that this process is parallel to applying
$\mbox{ls},i,
\mbox{lh},\delta$ to the RC
\[
\widetilde{RC}=(\mu,(\mu^{(2)},J^{(2)}),\ldots,(\mu^{(n)},J^{(n)})).
\]
Namely, if $b_l=\boxed{\scriptstyle{1}  \atop a}$, then the first letter to be obtained from
$\widetilde{RC}$ is
$a-1$. By doing these procedures successively until we obtain
$b_1,\ldots,b_l$, we can obtain
(3),(4). ($\bar{b}_1\ot\cdots\ot\bar{b}_l\in\P_+(\mu)$ follows from the fact
that it corresponds
to an admissible RC.)

These procedures continue until we finish removing letters corresponding to
$B^{(2)}_\mu$.
We then see $\tilde{p}$ corresponds to $(\mu^{(1)},J^{(1)})$ and therefore
has no letter
greater than 2. $\tilde{p}\in\P_+(\la)$ follows from the simple fact on
crystals that
if $b_1\ot b_2$ is a highest element, then so is $b_1$.
\end{proof}

Let $p$ be an element of $\P_+(\la)$. In view of the above theorem we
consider the following
procedure for $p$.
\begin{itemize}
\item[(1)] Supply $p$ with sufficiently many $u^{1,1}$ on the right.
\item[(2)] Compute the image of the map
\[
\begin{array}{ccc}
(B^{2,1})^{\ot L}\ot B_\la&\longrightarrow&B_\la\ot(B^{2,1})^{\ot L} \\
(u^{2,1})^{\ot L}\ot p&\longmapsto&\tilde{p}\ot p'.
\end{array}
\]
$L$ should be chosen so large that $\tilde{p}$ does not contain letters
greater than 2.
\item[(3)] Compute the RC $(\mu,J)$ corresponding to $\tilde{p}$.
\item[(4)] Cut the first rows of $p'$
(for their entries are all 1), decrease each letter by 1 and denote
it by $p''$. Compute
the image of the map
\[
\begin{array}{ccc}
B_\mu\ot(B^{1,1})^{\ot L}&\longrightarrow&(B^{1,1})^{\ot L}\ot B_\mu \\
(u^{1,\mu_1}\ot\cdots\ot u^{1,\mu_l})\ot p''&\longmapsto&u\ot p'''.
\end{array}
\]
(In fact, $u=(u^{1,1})^{\ot L}$.)
\item[(5)] Replace $p$ with $p'''$ and repeat (1)-(4) until $p'''$ has no
letter greater than 2.
\end{itemize}

This procedure is shown to give the RC corresponding to $p$ by 
Proposition \ref{prop:time evol},
Theorem 5.7 of \cite{Ta}, Theorem \ref{th:type2} and the
following lemma, that can easily be obtained from \cite{Sh}.

\begin{lemma}
Let $b\in B^{2,l},b'\in B^{2,m}$. Suppose the letters of the first rows of
$b,b'$ are all 1.
Let $\tilde{b}'\ot\tilde{b}$ be the image of $b\ot b'$ under the map
$B^{2,l}\ot B^{2,m}\rightarrow
B^{2,m}\ot B^{2,l}$. Then the letters of the first rows of
$\tilde{b}',\tilde{b}$ are again all 1,
and their second rows are given by the image of $b_2\ot b'_2$ under
$B^{1,l}\ot B^{1,m}\rightarrow
B^{1,m}\ot B^{1,l}$, where $b_2,b'_2$ are the second rows of $b,b'$.
\end{lemma}

\begin{example}\label{ex:okado}
$n=3$. We abbreviate $\otimes$ to $\cdot$ and 
draw frames of tableaux only for $B^{2,l}$.
\[
p=1\cd1\cd1\cd1\cd2\cd2\cd3\cd2\cd1\cd4\cd3\cd2\cd2\in B_{(1^{13})}.
\]
By $(B^{2,1})^{\ot3}\ot B_{(1^{16})}\rightarrow
B_{(1^{16})}\ot(B^{2,1})^{\ot3}$,
$(u^{2,1})^{\ot3}\ot(p\ot(u^{1,1})^{\ot3})$ maps to $\tilde{p}\ot p'$ where
\[
\tilde{p}=1\cd1\cd1\cd1\cd2\cd2\cd1\cd2\cd2\cd1\cd1\cd1\cd2\cd2\cd2\cd2,
\quad p'=\boxed{\scriptstyle{1} \atop3}\cd\boxed{\scriptstyle{1} \atop3}\cd
\boxed{\scriptstyle{1}  \atop 4}.
\]
The RC for $\tilde{p}$ is $(\mu,J)=((4,3,1),(0,1,4))$, and $p''=2\cd2\cd3$.
By
$B_\mu\ot(B^{1,1})^{\ot3}\rightarrow (B^{1,1})^{\ot3}\ot B_\mu$,
$(u^{1,4}\ot u^{1,3}\ot u^{1,1})\ot p''$ maps to $1^{\ot3}\ot p'''$ where
\[
p'''=1111\cd122\cd3.
\]
Redefine $p=p'''$ and repeat the procedure. We obtain
\begin{align}
&\tilde{p}=1111\cd122\cd1\cd2,\quad p'=
\boxed{\scriptstyle{1} \atop 3}, \\
&(\mu,J)=((2,1),(0,0)),\quad p''=2,\quad p'''=11\cd2.
\end{align}
The RC for this $p'''$ is $((1),(0))$. Thus we find the RC for the original
$p$ as
\begin{equation*}
{\small
\unitlength 10pt
\begin{picture}(20,6)(4,-2)
\multiput(6,2)(1,0){2}{\line(0,-1){3}}
\put(8,2){\line(0,-1){2}}
\put(9,2){\line(0,-1){2}}
\put(10,2){\line(0,-1){1}}

\put(6,2){\line(1,0){4}}
\put(6,1){\line(1,0){4}}
\put(6,0){\line(1,0){3}}
\put(6,-1){\line(1,0){1}}

\put(10.3,1.2){0}
\put(9.3,0.05){1}
\put(7.3,-0.95){4}
\put(7.5,3){$\mu^{(1)}$}
\multiput(13,0)(1,0){2}{\line(0,1){2}}
\put(15,2){\line(0,-1){1}}

\multiput(13,1)(0,1){2}{\line(1,0){2}}
\put(13,0){\line(1,0){1}}
\put(15.3,1.2){0}
\put(14.3,0.05){0}
\put(13.7,3){$\mu^{(2)}$}
\multiput(19,1)(1,0){2}{\line(0,1){1}}
\multiput(19,1)(0,1){2}{\line(1,0){1}}
\put(20.4,1.2){0}
\put(19,3){$\mu^{(3)}$}
\end{picture}
}
\end{equation*}

\end{example}

\section{Conjectures on other types}\label{sec:conj}

There is a bijection between 
the set of rigged configurations 
${\rm RC}(\lambda)$ and the highest paths (\ref{eq:Bl}) 
for $\geh_n = B^{(1)}_n$, $C^{(1)}_n$, $D^{(1)}_n$, $A^{(2)}_{2n-1}$,
$A^{(2)}_{2n}$ and $ D^{(2)}_{n+1}$.
It has been established in \cite{OSS, SS} by
a combinatorial algorithm 
similar to the $A^{(1)}_n$ case.
Here we state conjectures analogous to (\ref{eq:p}).

Precise specifications of the data ${\rm RC}(\lambda)$, 
${\mathcal C}_a, \Phi_a, p^{(n)}$ and $\bar{p}^{(n)}$ will be given 
in subsequent sections.
In particular for $\geh_n = A^{(2)}_{2n}$ and $D^{(2)}_{n+1}$,
we utilize the embedding $\iota$ into $C^{(1)}_n$ 
specified in (\ref{eq:embIII}) and (\ref{eq:iota}).
The relevant combinatorial $R$ for these algebras has been 
obtained in \cite{HKOT1,HKOT2} by an insertion scheme.
A piecewise linear formula is also available for 
$D^{(1)}_n$ in \cite{KOTY} and the other $\geh_n$ 
using the embedding described in \cite{KTT}.
\begin{conjecture}\label{conj:other}
Given a rigged configuration $(\mu, J)$ as in (\ref{eq:rc}), 
let $p \in {\mathcal P}_+(\lambda)$ 
be the classically highest path (\ref{eq:P}) that 
corresponds to it under the bijection 
in \cite{OSS,SS}. Then the following formulas are valid:
\begin{align}
p&=\Phi_1{\mathcal C}_1 \cdots \Phi_{n-2}{\mathcal C}_{n-2}
(\Phi_{n-1}{\mathcal C}_{n-1}(p^{(n)})
+ \Phi_n{\mathcal C}_n(\bar{p}^{(n)}))\;\;\hbox{ for } \;
\geh_n = D^{(1)}_n,\label{eq:pd}\\
p&= \Phi_1{\mathcal C}_1 \cdots 
\Phi_n{\mathcal C}_n(\bar{p}^{(n)})\;\; \qquad\qquad\qquad 
\qquad\hbox{for } \;
\geh_n = A^{(2)}_{2n-1}, B^{(1)}_n, C^{(1)}_n.\label{eq:pabc}
\end{align}
For $ \geh_n = A^{(2)}_{2n}$ and $D^{(2)}_{n+1}$,
let $p'$ be the above $\Phi_1{\mathcal C}_1 \cdots 
\Phi_n{\mathcal C}_n(\bar{p}^{(n)})$ for $C^{(1)}_n$ corresponding 
to the rigged configuration $\iota((\mu,J))$. Then 
$\iota^{-1}(p')$ exists and is equal to $p$. 
\end{conjecture}

For $\geh_n = A^{(2)}_{2n-1}, B^{(1)}_n$ and $C^{(1)}_n$,
define the $\geh_{n-a}$ path 
$p^{(a)} = \Phi_{a+1}{\mathcal C}_{a+1} \cdots 
\Phi_n{\mathcal C}_n(\bar{p}^{(n)})$ for 
$0 \le a \le n-1$.
Similarly for $\geh_n = D^{(1)}_n$, set  
$
p^{(a)} = \Phi_{a+1}{\mathcal C}_{a+1}\cdots 
\Phi_{n-2}{\mathcal C}_{n-2}$ $
(\Phi_{n-1}{\mathcal C}_{n-1}(p^{(n)})
+ \Phi_n{\mathcal C}_n(\bar{p}^{(n)}))
$
for $0 \le a \le n-2$.
Under Conjecture \ref{conj:other}, 
$p^{(a)}$ coincides with the image of the 
truncated rigged configuration 
$(\mu^{(a)}, (\mu^{(a+1)},J^{(a+1)}), \ldots, (\mu^{(n)},J^{(n)}))$ 
under the bijection in \cite{OSS,SS}.

\section{$D^{(1)}_n$ case}\label{sec:D(1)}

\subsection{Rigged configurations}\label{subsec:rcD(1)}

Consider the data of the form (\ref{eq:rc}) and 
define $E^{(a)}_j$ as in (\ref{eq:Eaj}).
The vacancy number is specified by (\ref{eq:paj}) for 
$1 \le a \le n-3$ and 
\begin{equation}
\begin{split}
p^{(n-2)}_j &= E^{(n-3)}_j+E^{(n-1)}_j+E^{(n)}_j -2 E^{(n-2)}_j,\\
p^{(a)}_j &= E^{(n-2)}_j -2 E^{(a)}_j\quad (a=n-1, n).
\end{split}
\end{equation}
The data (\ref{eq:rc}) is called a 
$D^{(1)}_n$ rigged configuration if (\ref{eq:cond}) is satisfied.
It is depicted as in $A^{(1)}_n$ case by the Young diagrams
with rigging.
The following is an example of $D^{(1)}_4$ rigged configuration.

\begin{equation}\label{eq:rcD(1)}
{\small
\unitlength 10pt
\begin{picture}(20,4)(4,-1)
\put(0,1){$(1^9)$}
\put(0,3){$\mu^{(0)}$}
\put(4,0){2}
\put(4,1){0}
\multiput(5,0)(1,0){4}{\line(0,1){2}}
\put(5,0){\line(1,0){3}}
\multiput(5,1)(0,1){2}{\line(1,0){4}}
\put(9,1){\line(0,1){1}}
\put(8.3,0.05){1}
\put(9.3,1.05){0}
\put(6,3){$\mu^{(1)}$}
\put(12,0){0}
\put(12,1){0}
\multiput(13,0)(1,0){3}{\line(0,1){2}}
\put(13,0){\line(1,0){2}}
\multiput(13,1)(0,1){2}{\line(1,0){3}}
\put(16,1){\line(0,1){1}}
\put(15.3,0.05){0}
\put(16.3,1.05){0}
\put(13.7,3){$\mu^{(2)}$}
\put(19,1){0}
\multiput(20,1)(1,0){3}{\line(0,1){1}}
\multiput(20,1)(0,1){2}{\line(1,0){2}}
\put(22.4,1.05){0}
\put(20.5,3){$\mu^{(3)}$}
\put(25,1){0}
\multiput(26,1)(1,0){3}{\line(0,1){1}}
\multiput(26,1)(0,1){2}{\line(1,0){2}}
\put(28.4,1.05){0}
\put(26.5,3){$\mu^{(4)}$}
\end{picture}
}
\end{equation}

\subsection{\mathversion{bold} Crystals}

\begin{equation}\label{eq:BD(1)}
B_l =\{ (x_1,\ldots,x_n,\bar{x}_n,\ldots, \bar{x}_1)
\in (\Z_{\ge 0})^{2n} \mid 
\sum_{i=1}^n(x_i+\bar{x}_i) = l, x_n\bar{x}_n =0 \}.
\end{equation}
\begin{equation}
B^{\ge a+1}_l = \{ (x_1,\ldots,x_n,\bar{x}_n,\ldots, \bar{x}_1)
 \in B_l \mid 
x_i=\bar{x}_i=0 \, (1 \le i \le a)\}
\end{equation}
for $0 \le a \le n-2$. Then 
\begin{equation*}
B_l = B^{\ge 1}_l \supset B^{\ge 2}_l \supset \cdots 
\supset B^{\ge n-1}_l 
\end{equation*}
as sets.
We regard $B^{\ge a+1}_l$ as $D^{(1)}_{n-a}$ crystal 
and represent the elements by 
length $l$ semistandard row tableaux over the alphabet 
$a+1 < \cdots <  n < \bar{n}<  \cdots < \overline{a+1}$.
(Actually the letter $n$ and $\bar{n}$ do not coexist in a 
tableau.)
We use the notation (\ref{eq:al}).
Then the highest element in $B^{\ge a+1}_l$ with respect to 
$D_{n-a}$ is $(a+1)^l$.
The rigged configuration (\ref{eq:rcD(1)}) corresponds to 
\begin{equation}\label{eq:D(1)path}
\unitlength 10pt
\begin{picture}(18.3,1.1)(1,0.1)
\multiput(-0.3,0.1)(2.3,0){9}{$\boxed{\phantom{3}}$}
\multiput(1.1,0.1)(2.3,0){8}{$\otimes$}
\put(0,0){1}
\put(2.3,0){$1$}
\put(4.6,0){$1$}
\put(6.9,0){$1$}
\put(9.2,0){$2$}
\put(11.5,0){$\bar{1}$}
\put(13.8,0){$\bar{1}$}
\put(16.1,0){$3$}
\put(18.4,0){$2\;\; $}
\end{picture}\, \in {\mathcal P}_+((1^9)).
\end{equation}

\subsection{\mathversion{bold} 
${\mathcal C}_{n-1}(p^{(n)})$ and 
${\mathcal C}_{n}(\bar{p}^{(n)})$}
\label{subsec:cpcp}

Set
\begin{align}
p^{(n)}
&= n^{\mu^{(n-1)}_1}\otimes \cdots \otimes 
n^{\mu^{(n-1)}_{l_{n-1}}},\\
\bar{p}^{(n)}&
= \bar{n}^{\mu^{(n)}_1}\otimes \cdots \otimes 
\bar{n}^{\mu^{(n)}_{l_{n}}},\label{eq:pnbar}
\end{align}
where $n^l $ means the unique element 
of ``$A^{(1)}_0$ crystal" $B_l$ realized by the tableau
with letter $n$ only.  (See (\ref{eq:al}).) 
$\bar{n}^l$ stands for the same object in another 
copy of $A^{(1)}_0$ crystal $B_l$.
Then 
${\mathcal C}_{n-1}(p^{(n)})$ is defined from 
(\ref{eq:Ca}) and (\ref{eq:di}) by setting  
$(\mu, J) = (\mu^{(n-1)}, J^{(n-1)})$, 
$b_i = n^{\mu_i} (i \ge 1)$, $b_0 = n^{\mu_1}$ and 
regarding $B_l$ there as the $A^{(1)}_0$ crystal with letter $n$ only.
Namely, 
${\mathcal C}_{n-1}(p^{(n)}) = :\!n^{\mu_1}[d_1]
\otimes \cdots \otimes
n^{\mu_m}[d_m]\!:$, where the normal ordering $::$ is done 
along the $A^{(1)}_0$ combinatorial $R$ 
(\ref{eq:A0}) with $n+1$ replaced by $n$.
Similarly ${\mathcal C}_{n}(\bar{p}^{(n)})$
is defined by setting  
$(\mu, J) = (\mu^{(n)}, J^{(n)})$, 
$b_i = {\bar n}^{\mu_i} (i \ge 1)$ and $b_0 = \bar{n}^{\mu_1}$
in (\ref{eq:Ca}) and (\ref{eq:di}).
As the result 
${\mathcal C}_{n-1}(p^{(n)})$ and ${\mathcal C}_{n}(\bar{p}^{(n)})$
are elements of the two independent copies of 
$A^{(1)}_0$ affine crystals realized with tableaux with letters 
$n$ only or $\bar{n}$ only.
Since the energy function in (\ref{eq:A0}) is positive,
the $i$-th mode $d_i$ in the normal ordered form satisfies
$0 \le d_1 \le \cdots \le d_m$.

\subsection{\mathversion{bold} 
$\Phi_{n-1}{\mathcal C}_{n-1}(p^{(n)})$,  
$\Phi_n{\mathcal C}_{n}(\bar{p}^{(n)})$ 
and their superposition.}\label{subsec:super}

Let us write as 
$(\lambda_1,\ldots,\lambda_k)
=(\mu^{(n-2)}_1,\ldots, \mu^{(n-2)}_{l_{n-2}})$
and 
${\mathcal C}_{n-1}(p^{(n)}) 
= n^{\mu_1}[d_1]\otimes \cdots \otimes n^{\mu_m}[d_m]$.
Let $B'_l=\{(x_{n-1},x_n)' \in (\Z_{\ge 0})^2\mid x_{n-1}+x_n=l\}$ 
be the $A^{(1)}_1$ crystal with tableau letters $n-1,n$,  
where $n^l=\boxed{n\cdots n}$ 
is identified as the lowest weight element
$(0,l)'$.
We define $\Phi_{n-1}{\mathcal C}_{n-1}(p^{(n)})
\in B'_{\lambda_1}\otimes \cdots \otimes B'_{\lambda_k}$ 
to be $c_1\otimes \cdots \otimes c_k$ determined 
by (\ref{eq:zdon}) by setting $b_i = n^{\mu_i}$ and $a=n-1$.

$\Phi_{n}{\mathcal C}_{n}(\bar{p}^{(n)})$ is defined similarly by 
replacing the letter $n$ by $\bar{n}$.
To be precise, we reset the meaning of $\mu_i$ and $d_i$
by saying that  ${\mathcal C}_{n}(\bar{p}^{(n)}) 
= {\bar n}^{\mu_1}[d_1]\otimes \cdots \otimes \bar{n}^{\mu_m}[d_m]$.
Let $B''_l=\{(x_{n-1},\bar{x}_n)'' 
\in (\Z_{\ge 0})^2\mid x_{n-1}+\bar{x}_n=l\}$ 
be the $A^{(1)}_1$ crystal with tableau letters $n-1,\bar{n}$,  
where $\bar{n}^l = 
\begin{picture}(35,1)
\put(0,0){$\boxed{\phantom{n\cdots n}}$}
\put(3,0){$\bar{n}\cdots\bar{n}$}
\end{picture}
$
is identified as the lowest weight element $(0,l)''$.
We keep the meaning $\lambda=\mu^{(n-2)}$ as before.
Then $\Phi_{n}{\mathcal C}_{n}(\bar{p}^{(n)}) 
\in B''_{\lambda_1}\otimes \cdots \otimes B''_{\lambda_k}$ is defined to be
$c_1\otimes \cdots \otimes c_k$ determined 
by (\ref{eq:zdon}) by setting $b_i = \bar{n}^{\mu_i}$ and $a=n-1$.

We introduce the ``superposition" $+$ of the 
two copies of $A^{(1)}_1$ crystals $B'_l$ and $B''_l$:
\begin{equation}\label{eq:super}
\begin{split}
& \qquad\qquad 
\qquad +: \; B'_l  \times B''_l \qquad\qquad \longrightarrow 
\qquad D^{(1)}_{2} \; \hbox{ crystal} \;\; B^{\ge n-1}_l \\
&\qquad\qquad \;(x_{n-1},x_n)' \times (y_{n-1},\bar{y}_n)'' 
\; \;\longmapsto \;
(0, \ldots, 0, z_{n-1}, z_n, \bar{z}_n, \bar{z}_{n-1}, 0, \ldots, 0),\\
&z_{n-1} = \min(x_{n-1}, y_{n-1}), 
z_n = \max(0, y_{n-1}-x_{n-1}),\\
&\bar{z}_n =\max(0, \bar{y}_n - x_n), 
\bar{z}_{n-1} = \min(x_n, \bar{y}_n).
\end{split}
\end{equation}
In terms of the tableaux, this $+$ means to ``superpose" the 
letters $n-1,n$ from $B'_l$ and $n-1, \bar{n}$ from $B''_l$ 
as $(n-1,n-1) \mapsto n-1,
(n,n-1) \mapsto n,
(n-1,\bar{n}) \mapsto \bar{n}$ and 
$(n,\bar{n}) \mapsto \overline{n-1}$.
For example in $n=4$ case, one has 
$\boxed{334}+
\begin{picture}(23,1)
\unitlength 10pt 
\put(0,0){$\boxed{\phantom{444}}$}
\put(0.3,0){$3\bar{4}\bar{4}$}
\end{picture}
= 
\begin{picture}(23,1)
\unitlength 10pt 
\put(0,0){$\boxed{\phantom{444}}$}
\put(0.3,0){$3\bar{4}\bar{3}$}
\end{picture}
$.
Note that the conditions $z_n\bar{z}_n=0$ and 
$z_{n-1}+z_n+\bar{z}_n+\bar{z}_{n-1}=l$ are satisfied. 
Under these conditions, the map (\ref{eq:super}) 
is invertible.
Suppose we have 
\begin{align*}
\Phi_{n-1}{\mathcal C}_{n-1}(p^{(n)})  &
= b'_1\otimes \cdots \otimes b'_k 
\in B'_{\lambda_1}\otimes \cdots \otimes B'_{\lambda_k},\\
\Phi_{n}{\mathcal C}_{n}(\bar{p}^{(n)}) &
= b''_1\otimes \cdots \otimes b''_k 
\in B''_{\lambda_1}\otimes \cdots \otimes B''_{\lambda_k},
\end{align*}
where $\lambda = \mu^{(n-2)}$ as before.
Then we specify their superposition as
\begin{equation}\label{eq:plus}
\Phi_{n-1}{\mathcal C}_{n-1}(p^{(n)}) +
\Phi_{n}{\mathcal C}_{n}(\bar{p}^{(n)}) =
(b'_1+b''_1)\otimes \cdots \otimes (b'_k + b''_k)
\in B^{\ge n-1}_{\lambda_1}\otimes \cdots \otimes 
B^{\ge n-1}_{\lambda_k}.
\end{equation}

\subsection{\mathversion{bold} Maps 
${\mathcal C}_1, \ldots, {\mathcal C}_{n-2}$}

For $1 \le a \le n-2$, ${\mathcal C}_a$ are 
given  by (\ref{eq:Ca}) and (\ref{eq:di}) 
provided that $B_l$ is regarded as the 
$D^{(1)}_{n-a}$ crystal $B^{\ge a+1}_l$.
The normal ordered form 
$b[d_1]\otimes \cdots \otimes b_m[d_m]$ 
is defined in the same way 
as the $A^{(1)}_n$ case described in 
Section \ref{subsec:noA(1)}.
A new feature is that the energy function can be negative,
hence $d_1 \le \cdots \le d_m$ is no longer valid in general.

\subsection{\mathversion{bold} Maps 
$\Phi_1, \ldots, \Phi_{n-2}$}\label{subsec:PhiD}

For $1\le a \le n-2$, we define the map 
$\Phi_a$ by (\ref{eq:Phi}) and (\ref{eq:zdon}) 
by regarding $B_l = B^{\ge a+1}_l$ and 
$B'_l = B^{\ge a}_l$ as the $D^{(1)}_{n-a}$ and 
$D^{(1)}_{n-a+1}$ crystals, provided 
$d_1 \le \cdots \le d_m$ holds
among the modes of the normal ordered element
$b_1[d_1]\ot \cdots \ot b_m[d_m]$.
In contrast to the $A^{(1)}_n$ case, 
this is no longer valid in general 
since the energy function can be negative.
($d_1 \ge 0$ can be assured.)
To make sense of ${\mathcal T}^{d}_a=\boxed{a}^{\otimes d}$ for  
negative $d$ in (\ref{eq:zdon}), 
we make the regularization explained in the sequel.

In (\ref{eq:zdon}), take the largest $i$ such that 
$d_i > d_{i+1}$. Then its left hand side looks as
\begin{equation}\label{eq:pre0}
({\mathcal T}_a^{d_1}\otimes
\cdots  \otimes b_i \otimes 
{\mathcal T}_a^{-\delta} \ot b_{i+1} \ot 
{\mathcal T}_a^{d_{i+2}-d_{i+1}}\otimes \cdots \ot b_m)\ot 
(a^{\lambda_1} \ot a^{\lambda_2} \ot \cdots \ot a^{\lambda_k}),
\end{equation}
which is regarded as an element in the 
tensor product of $D^{(1)}_{n-a+1}$ crystal $B^{\ge a}_l$'s.
We have set $\delta= d_i-d_{i+1}>0$. 
There is no problem in carrying the components 
$(b_{i+1} \otimes \cdots \otimes b_m)$ through 
$(a^{\lambda_1} \otimes \cdots \otimes a^{\lambda_k})$
to the right to find $p_j$'s in 
$(b_{i+1} \otimes \cdots \otimes b_m)\otimes 
(a^{\lambda_1} \otimes \cdots \otimes a^{\lambda_k})
\simeq (p_1 \otimes \cdots \otimes p_k) \otimes {\rm tail}$.
Thus we are to define 
${\mathcal T}^{-\delta}_a \otimes 
(p_1 \otimes \cdots \otimes p_k)\otimes {\rm tail}$.
Actually we declare that 
$b_i \otimes {\mathcal T}^{-\delta}_a \otimes 
(p_1 \otimes \cdots \otimes p_k)\otimes {\rm tail}$ 
is to be understood as
\begin{equation}\label{eq:pre1}
{\tilde b}_i \otimes 
({\tilde p}_1 \otimes \cdots \otimes {\tilde p}_k) 
\otimes {\rm tail},
\end{equation}
where ${\tilde p}_1 \otimes \cdots \otimes {\tilde p}_k$ is 
determined by sending ${\mathcal T}^\delta_a$ from 
the right to the left by the $D^{(1)}_{n-a+1}$ combinatorial $R$:
\begin{equation}
(p_1 \otimes \cdots \otimes p_k) \otimes \boxed{a}^{\,\otimes \delta}
\simeq {\rm head} \otimes 
({\tilde p}_1 \otimes \cdots \otimes {\tilde p}_k).
\end{equation}

\vspace{0.2cm}\noindent
In the actual use, ${\rm head} = 
\begin{picture}(35,1.5)
\put(0,0){$\boxed{\phantom{a\!+\!1}}^{\,\otimes \delta}$}
\put(3,-0.2){$\overline{a\!+\!1}$}
\end{picture}$ 
seems valid, but this is not necessary for our definition.
To specify ${\tilde b}_i$ in (\ref{eq:pre1}),  we use the expression 
$b_i= (x_1,\ldots, x_n,{\bar x}_n,\ldots  {\bar x}_1)
\in B^{\ge a+1}_{\mu_i}$.
Then ${\tilde b}_i$ is given by modifying it as 
$(x_{a+1},{\bar x}_a=0) \rightarrow (x_{a+1}-\delta, \delta)$ leaving 
the other coordinates unchanged.
It is a part of our conjecture that $x_{a+1} \ge \delta$ is always valid 
as the result of normal ordering hence ${\tilde b}_i$ is well defined.

This regularization of $b_i \otimes {\mathcal T}_a^{d_{i+1}-d_i}$ 
with $d_i> d_{i+1}$ can be done successively also in the 
part $({\mathcal T}^{d_1}_a \otimes b_1 \otimes \cdots 
\otimes b_{i-1} \otimes {\mathcal T}_a^{d_i-d_{i-1}})$
in (\ref{eq:pre0}), which leads to the (conjectural) definition of 
our $\Phi_a$.
The change $b_i \rightarrow {\tilde b}_i$ is analogous 
to the pair annihilation and creation in the 
$D^{(1)}_n$ automaton described in \cite{KTT}.

\subsection{Example}\label{subsec:exD(1)}

Consider the rigged configuration (\ref{eq:rcD(1)}).
We have $p^{(4)}=\boxed{44}$, 
$\bar{p}^{(4)} = 
\begin{picture}(18,1)
\put(0,0){$\boxed{\phantom{44}}$}
\put(3,0){$\bar{4}\bar{4}$}
\end{picture}
$, 
${\mathcal C}_3(p^{(4)})=\boxed{44}_{\,2}$ and 
${\mathcal C}_4(\bar{p}^{(4)}) = 
\begin{picture}(23,1)
\put(0,0){$\boxed{\phantom{44}}_{\,2}$}
\put(3,0){$\bar{4}\bar{4}$}
\end{picture}
$.
{}From 
\begin{equation*}
\unitlength 1mm
\begin{picture}(30,40)(0,-7)
\multiput(0,0)(0,12){3}{
\multiput(0,0)(16,0){2}{\line(1,0){10}}
\multiput(5,-3)(16,0){2}{\line(0,1){6}}
}

\put(2,29){333}
\put(19,29){33}

\put(-5,23){44}
\put(11,23){33}
\put(27,23){33}

\put(2,17){344}
\put(19,17){33}

\put(-4,11){3}
\put(12,11){4}
\put(28,11){3}

\put(2,5){334}
\put(19,5){34}

\put(-4,-1){3}
\put(12,-1){4}
\put(28,-1){3}

\put(2,-6){333}
\put(19,-6){44}

\end{picture}
\end{equation*}
we have $\Phi_3{\mathcal C}_3(p^{(4)}) = \boxed{333}\otimes \boxed{44}$.
Similarly, 
$\Phi_4{\mathcal C}_4(\bar{p}^{(4)}) = 
\boxed{333}\otimes 
\begin{picture}(18,1)
\unitlength 10pt 
\put(0,0){$\boxed{\phantom{44}}$}
\put(0.3,0){$\bar{4}\bar{4}$}
\end{picture}
$.
Thus the superposition is determined as 
$p^{(2)}= \Phi_3{\mathcal C}_3(p^{(4)})+
\Phi_4{\mathcal C}_4(\bar{p}^{(4)})=
\boxed{333}\otimes
\begin{picture}(18,1)
\unitlength 10pt 
\put(0,0){$\boxed{\phantom{44}}$}
\put(0.3,0){$\bar{3}\bar{3}$}
\end{picture}
$.
Since 
\begin{equation*}
\boxed{333}_{\,3}\otimes
\begin{picture}(18,1)
\unitlength 10pt 
\put(0,0){$\boxed{\phantom{44}}_{\,0}$}
\put(0.3,0){$\bar{3}\bar{3}$}
\end{picture}
\;
\simeq 
\boxed{33}_{\,2} \otimes 
\begin{picture}(28,1)
\unitlength 10pt 
\put(0,0){$\boxed{\phantom{333}}_{\,1}$}
\put(0.3,0){$3\bar{3}\bar{3}$}
\end{picture}
\end{equation*}
in $D^{(1)}_2$ crystal, we obtain ${\mathcal C}_2(p^{(2)}) = 
\boxed{33}_{\,2}\otimes
\begin{picture}(28,1)
\unitlength 10pt 
\put(0,0){$\boxed{\phantom{333}}_{\,1}$}
\put(0.3,0){$3\bar{3}\bar{3}$}
\end{picture}
$.
Next to find $p^{(1)}=\Phi_2{\mathcal C}_2(p^{(2)})$,
we compute 
${\mathcal T}_2^2\otimes \boxed{33}\otimes
{\mathcal T}_2^{-1} \otimes 
\begin{picture}(23,1)
\unitlength 10pt 
\put(0,0){$\boxed{\phantom{333}}$}
\put(0.3,0){$3\bar{3}\bar{3}$}
\end{picture} \otimes (\boxed{2222}\otimes \boxed{222})
$ in $D^{(1)}_3$ crystal as follows:
\begin{equation*}
\unitlength 1mm
\begin{picture}(30,65)(0,-32)
\multiput(0,-24)(0,12){5}{
\multiput(0,0)(18,0){2}{\line(1,0){10}}
\multiput(5,-3)(18,0){2}{\line(0,1){6}}
}

\put(-17,-1){33}\put(-11.5,-1){$\rightsquigarrow$}
\put(-5,-1){$3\bar{2}$}
\put(12,-1){23}
\put(29,-1){22}

\put(-4,11){$\bar{3}$}
\put(13,11){2}
\put(30,11){2}\put(25,11.2){$<$}

\put(-5,12){\line(-1,0){5}}\put(-10,12){\vector(0,-1){10}}

\put(-6,23){$3\bar{3}\bar{3}$}
\put(12,23){222}
\put(29,23){222}

\put(2,29){2222}
\put(20.5,29){222}

\put(2,17){$23\bar{3}\bar{3}$}
\put(20.5,17){222}

\put(2,5){$223\bar{3}$}
\put(20.5,5){222}

\put(2,-7){$22\bar{2}\bar{2}$}
\put(20.5,-7){223}

\put(-4,-13){2}
\put(13,-13){$\bar{2}$}
\put(30,-13){2}

\put(2,-19){$222\bar{2}$}
\put(20.5,-19){$23\bar{2}$}

\put(-4,-25){2}
\put(13,-25){$\bar{2}$}
\put(30,-25){2}

\put(2,-31){2222}
\put(20.5,-31){$3\bar{2}\bar{2}$}

\end{picture}
\end{equation*}
Here the symbol $<$ in the horizontal line
signifies ${\mathcal T}_2^{-1}$ where the 
combinatorial $R$ acts from NE to SW.
(Our usual convention is from NW to SE as in (\ref{eq:Rfig}). )
The change $33$ into $3\bar{2}$ due to the pair 
annihilation and creation is indicated by $\rightsquigarrow$.
We have $p^{(1)}=\boxed{2222}\otimes
\begin{picture}(23,1)
\unitlength 10pt 
\put(0,0){$\boxed{\phantom{333}}$}
\put(0.3,0){$32\bar{2}$}
\end{picture}$.
Under the $D^{(1)}_3$ combinatorial $R$, one has  
\begin{equation*}
\boxed{2222}_{\,4}\otimes
\begin{picture}(23,1)
\unitlength 10pt 
\put(0,0){$\boxed{\phantom{333}}_{\,2}$}
\put(0.3,0){$3\bar{2}\bar{2}$}
\end{picture}
\;
\simeq 
\boxed{222}_{\,4} \otimes 
\begin{picture}(28,1)
\unitlength 10pt 
\put(0,0){$\boxed{\phantom{3333}}_{\,2}\, .$}
\put(0.3,0){$23\bar{2}\bar{2}$}
\end{picture}
\end{equation*}

Therefore the both sides are normal ordered forms and 
we choose the left hand side as 
${\mathcal C}_1(p^{(1)})$.
Finally to obtain 
$p= \Phi_1{\mathcal C}_1(p^{(1)})$, 
we compute 
${\mathcal T}_1^4 \otimes \boxed{2222}  
\otimes {\mathcal T}_1^{-2} \otimes
\begin{picture}(25,1)
\unitlength 10pt 
\put(0,0){$\boxed{\phantom{333}}$}
\put(0.3,0){$3\bar{2}\bar{2}$}
\end{picture}
\otimes (\,\boxed{1}^{\otimes 9}\,)$ in 
$D^{(1)}_4$ crystal.
It is calculated by applying ${\mathcal T}_1^4\otimes$
to the output of the following diagram
($\,\boxed{1}^{\otimes 9}$ has been 
replaced with $\,\boxed{1}^{\otimes 6}$ below): 
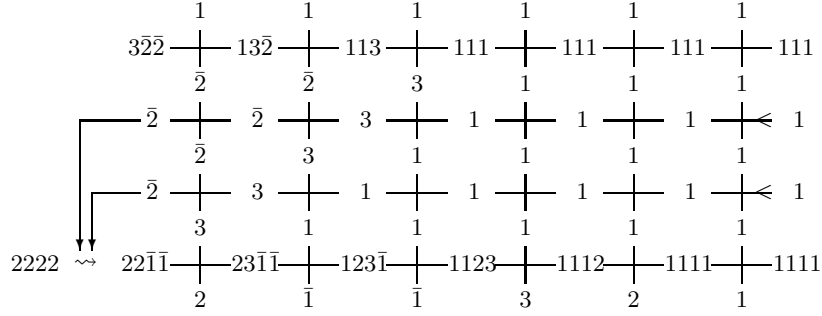
\begin{equation*}
\unitlength 0.8mm
{\small
\begin{picture}(110,53)(-10,-19)
\multiput(0,-12)(0,12){4}{
\multiput(0,0)(18,0){6}{\line(1,0){10}}
\multiput(5,-3)(18,0){6}{\line(0,1){6}}
}

\multiput(4,29)(18,0){6}{\put(0,0){1}}

\put(-7,23){$3\bar{2}\bar{2}$}
\put(11,23){$13\bar{2}$}
\put(29,23){113}
\multiput(47,23)(18,0){4}{\put(0,0){111}}

\put(4,17){$\bar{2}$}\put(22,17){$\bar{2}$}\put(40,17){3}
\multiput(58,17)(18,0){3}{\put(0,0){1}}

\put(-5,12){\line(-1,0){10}}\put(-15,12){\vector(0,-1){21.5}}
\put(-4,11){$\bar{2}$}\put(13.5,11){$\bar{2}$}\put(31.5,11){3}
\multiput(49.5,11)(18,0){4}{\put(0,0){1}}
\put(97,11.2){$<$}

\put(4,5){$\bar{2}$}\put(22,5){3}
\multiput(40,5)(18,0){4}{\put(0,0){1}}

\put(-5,0){\line(-1,0){8}}\put(-13,0){\vector(0,-1){9.5}}
\put(-4,-1){$\bar{2}$}\put(13.5,-1){3}\put(31.5,-1){1}
\multiput(49.5,-1)(18,0){4}{\put(0,0){1}}
\put(97,-0.8){$<$}

\put(4,-7){3}
\multiput(22,-7)(18,0){5}{\put(0,0){1}}

\put(-26.5,-13){2222}\put(-16,-12.5){$\rightsquigarrow$}
\put(-8.2,-13){$22\bar{1}\bar{1}$}
\put(10.2,-13){$23\bar{1}\bar{1}$}
\put(28.2,-13){$123\bar{1}$}
\put(46.2,-13){1123}
\put(64.2,-13){1112}
\put(82.2,-13){1111}
\put(100.2,-13){1111}

\put(4,-19){2}\put(22,-19){$\bar{1}$}\put(40,-19){$\bar{1}$}
\put(58,-19){3}\put(76,-19){2}\put(94,-19){1}

\end{picture}
}
\end{equation*}
Thus we get $p=11112\bar{1}\bar{1}32$ 
in agreement with (\ref{eq:D(1)path}).

\section{$A^{(2)}_{2n-1}$ case}\label{sec:A(2)odd}

The $A^{(2)}_{2n-1}$ case can be reduced to 
$D^{(1)}_{n+1}$ by dropping
one term in the superposition in (\ref{eq:pd}).

\subsection{Rigged configurations}

For the data (\ref{eq:rc}) define $E^{(a)}_j$ as in (\ref{eq:Eaj}).
The vacancy numbers are specified by (\ref{eq:paj}) 
for $1 \le a \le n-2$ and 
\begin{equation}
\begin{split}
p^{(n-1)}_j &=E^{(n-2)}_j-2E^{(n-1)}_j+2E^{(n)}_j,\\
p^{(n)}_j &=E^{(n-1)}_j-2E^{(n)}_j.
\end{split}
\end{equation}
The data (\ref{eq:rc}) is 
called an $A^{(2)}_{2n-1}$ rigged configuration
if the condition (\ref{eq:cond}) is satisfied.
The following is an example of $A^{(2)}_5$ rigged configuration.

\begin{equation}\label{eq:rcA(2)odd}
{\small
\unitlength 10pt
\begin{picture}(20,4)(4,-1)
\put(0,1){$(1^7)$}
\put(0,3){$\mu^{(0)}$}
\put(5,-0.3){3}
\put(5,1.2){3}
\multiput(6,2)(1,0){2}{\line(0,-1){3}}
\put(8,2){\line(0,-1){1}}
\put(6,2){\line(1,0){2}}
\put(6,1){\line(1,0){2}}
\multiput(6,0)(0,-1){2}{\line(1,0){1}}
\put(8.3,1.2){3}
\put(7.3,0.05){1}
\put(7.3,-0.95){3}
\put(6.5,3){$\mu^{(1)}$}
\put(12,0.7){0}
\multiput(13,0)(1,0){3}{\line(0,1){2}}
\put(13,0){\line(1,0){2}}
\multiput(13,1)(0,1){2}{\line(1,0){2}}
\put(15.3,1.2){0}
\put(15.3,0.05){0}
\put(13.7,3){$\mu^{(2)}$}
\put(19,1.2){0}
\multiput(20,1)(1,0){3}{\line(0,1){1}}
\multiput(20,1)(0,1){2}{\line(1,0){2}}
\put(22.4,1.2){0}
\put(20.5,3){$\mu^{(3)}$}
\end{picture}
}
\end{equation}

\subsection{Crystals}

\begin{align*}
B_l &=\{ (x_1,\ldots,x_n,\bar{x}_n,\ldots, \bar{x}_1)
\in (\Z_{\ge 0})^{2n} \mid 
\sum_{i=1}^n(x_i+\bar{x}_i) = l\},\\
B^{\ge a+1}_l &= \{ (x_1,\ldots,x_n,\bar{x}_n,\ldots, \bar{x}_1)
 \in B_l \mid 
x_i=\bar{x}_i=0 \, (1 \le i \le a)\}\quad (0 \le a \le n-1).
\end{align*}
\begin{equation*}
B_l = B^{\ge 1}_l \supset B^{\ge 2}_l \supset \cdots 
\supset B^{\ge n}_l 
\end{equation*}
as sets.
We regard $B^{\ge a+1}_l$ as $A^{(2)}_{2n-2a+1}$ crystal 
and represent the elements by 
length $l$ semistandard row tableaux over the alphabet 
$a+1 < \cdots <  n < \bar{n}<  \cdots < \overline{a+1}$.
We use the notation (\ref{eq:al}).
Then the highest element in $B^{\ge a+1}_l$ with respect to 
the classical part $C_{n-a}$ is $(a+1)^l$.
The rigged configuration (\ref{eq:rcA(2)odd}) 
corresponds to the following element in 
${\mathcal P}_+((1^7))$:
\begin{equation}\label{eq:pA(2)odd}
\unitlength 10pt
\begin{picture}(20,1.1)(-3,0.1)
\multiput(-0.3,0.1)(2.3,0){7}{$\boxed{\phantom{3}}$}
\multiput(1.1,0.1)(2.3,0){6}{$\otimes$}
\put(0,0){1}
\put(2.3,0){$1$}
\put(4.6,0){$2$}
\put(6.9,0){$1$}
\put(9.2,0){$1$}
\put(11.5,0){$\bar{2}$}
\put(13.8,0){$\bar{1}\;\;.$}
\end{picture}
\end{equation}

\subsection{\mathversion{bold} Maps 
${\mathcal C}_1, \ldots, {\mathcal C}_{n}$}

For $1 \le a \le n-1$, ${\mathcal C}_a$ is  
defined by (\ref{eq:Ca}) and (\ref{eq:di}) 
provided that $B_l$ is regarded as the 
$A^{(2)}_{2n-2a-1}$ crystal $B^{\ge a+1}_l$.

We define $\bar{p}^{(n)}$ and 
${\mathcal C}_n(\bar{p}^{(n)})$ in the same way as
in Section \ref{subsec:cpcp}.
The latter is obtained by assigning the modes 
to the former and 
normal ordering by using the 
$A^{(1)}_0$ combinatorial $R$ (\ref{eq:A0}) 
with letter $\bar{n}$.

\subsection{\mathversion{bold} Maps 
$\Phi_1, \ldots, \Phi_{n}$}

For $1\le a \le n-1$, we define the map 
$\Phi_a$ by (\ref{eq:Phi}) and (\ref{eq:zdon}) 
by regarding $B_l = B^{\ge a+1}_l$ and 
$B'_l = B^{\ge a}_l$ as the $A^{(2)}_{2n-2a-1}$ and 
$A^{(2)}_{2n-2a+1}$ crystals.
Since the energy function can be negative we employ
the same regularization as $D^{(1)}_n$ explained in 
Section \ref{subsec:PhiD}.

$\Phi_n{\mathcal C}_n(\bar{p}^{(n)})$ 
is defined similarly to Section \ref{subsec:super}
by changing the letter $n-1$ to $n$.
To be precise, let 
${\mathcal C}_{n}(\bar{p}^{(n)}) 
= {\bar n}^{\mu_1}[d_1]\otimes \cdots 
\otimes \bar{n}^{\mu_m}[d_m]$.
Regard $B^{\ge n}_l
=\{(0,\ldots,0,x_{n},\bar{x}_n,0,\ldots,0)
\in (\Z_{\ge 0})^2\mid x_{n}+\bar{x}_n=l\}$ 
as the $A^{(1)}_1$ crystal with tableau letters 
$n,\bar{n}$,  where $\bar{n}^l = 
\begin{picture}(35,1)
\put(0,0){$\boxed{\phantom{n\cdots n}}$}
\put(3,0){$\bar{n}\cdots\bar{n}$}
\end{picture}
$
is identified as the lowest weight element.
Set $\lambda=\mu^{(n-1)}$. Then 
$\Phi_{n}{\mathcal C}_{n}(\bar{p}^{(n)}) 
\in B^{\ge n}_{\lambda_1}\otimes 
\cdots \otimes B^{\ge n}_{\lambda_k}$ is defined to be
$c_1\otimes \cdots \otimes c_k$ determined 
by (\ref{eq:zdon}) with 
$b_i = \bar{n}^{\mu_i}$ and $a=n$.

\subsection{Example}

Consider the rigged configuration (\ref{eq:rcA(2)odd}).
We have $\bar{p}^{(3)} = 
\begin{picture}(18,1)
\put(0,0){$\boxed{\phantom{33}}$}
\put(3,0){$\bar{3}\bar{3}$}
\end{picture}$
and 
${\mathcal C}_3(\bar{p}^{(3)}) =
\begin{picture}(25,1)
\put(0,0){$\boxed{\phantom{33}}_{\, 2}$}
\put(3,0){$\bar{3}\bar{3}$}
\end{picture}$.
To find 
$\Phi_3{\mathcal C}_3(\bar{p}^{(3)})$, we calculate 
${\mathcal T}_3^2 \otimes 
\begin{picture}(18,1)
\put(0,0){$\boxed{\phantom{33}}$}
\put(3,0){$\bar{3}\bar{3}$}
\end{picture}
\otimes (\boxed{33} \otimes \boxed{33})$ as
\begin{equation*}
\unitlength 1mm
\begin{picture}(30,40)(0,-7)
\multiput(0,0)(0,12){3}{
\multiput(0,0)(16,0){2}{\line(1,0){10}}
\multiput(5,-3)(16,0){2}{\line(0,1){6}}
}

\put(3,29){33}
\put(19,29){33}

\put(-5,23){$\bar{3}\bar{3}$}
\put(11,23){33}
\put(27,23){33}

\put(3,17){$\bar{3}\bar{3}$}
\put(19,17){33}

\put(-4,11){3}
\put(12,11){$\bar{3}$}
\put(28,11){3}

\put(3,5){$3\bar{3}$}
\put(19,5){$3\bar{3}$}

\put(-4,-1){3}
\put(12,-1){$\bar{3}$}
\put(28,-1){3}

\put(3,-6){33}
\put(19,-6){$\bar{3}\bar{3}$}

\end{picture}
\end{equation*}
Thus we get $p^{(2)}=\Phi_3{\mathcal C}_3(\bar{p}^{(3)})
= \boxed{33}\otimes
\begin{picture}(18,1)
\put(0,0){$\boxed{\phantom{33}}$}
\put(3,0){$\bar{3}\bar{3}$}
\end{picture}
$.
In view of $H(\,\boxed{33}\otimes
\begin{picture}(18,1)
\put(0,0){$\boxed{\phantom{33}}\,)$}
\put(3,0){$\bar{3}\bar{3}$}
\end{picture} \,= -2$,
we have
${\mathcal C}_2({p}^{(2)}) = 
\boxed{33}_{\,2}\otimes
\begin{picture}(25,1)
\put(0,0){$\boxed{\phantom{33}}_{\, 0}$}
\put(3,0){$\bar{3}\bar{3}$}
\end{picture}
$.
To find $\Phi_2{\mathcal C}_2(p^{(2)})$, 
we calculate 
${\mathcal T}_2^2\otimes \boxed{33}
\otimes {\mathcal T}^{-2}_2 \otimes 
\begin{picture}(18,1)
\put(0,0){$\boxed{\phantom{33}}$}
\put(3,0){$\bar{3}\bar{3}$}
\end{picture}
\otimes (\boxed{22} \otimes \boxed{2} \otimes \boxed{2})
$
by using the $A^{(2)}_3$ combinatorial $R$ 
on letters $2,3,\bar{3},\bar{2}$ as
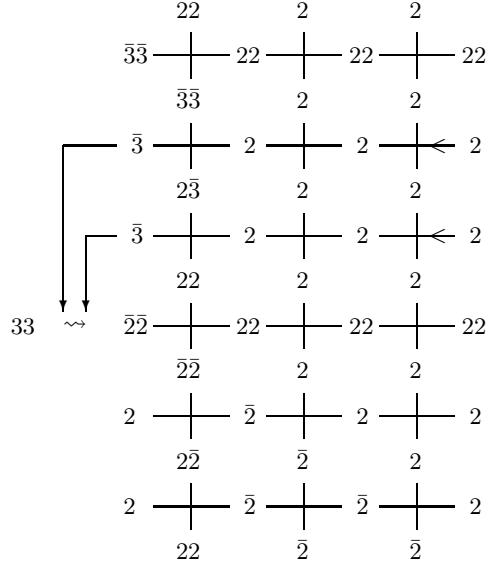
\begin{equation*}
{\small 
\unitlength 1mm
\begin{picture}(30,77)(5,-32)
\multiput(0,-24)(0,12){6}{
\multiput(0,0)(15,0){3}{\line(1,0){10}}
\multiput(5,-3)(15,0){3}{\line(0,1){6}}
}

\put(3,41){22}
\put(19,41){2}
\put(34,41){2}

\put(-4,35){$\bar{3}\bar{3}$}
\put(11,35){$22$}
\put(26,35){22}
\put(41,35){22}

\put(3,29){$\bar{3}\bar{3}$}
\put(19,29){2}
\put(34,29){2}

\put(-3,23){$\bar{3}$}\put(-5,24){\line(-1,0){7}}
\put(12,23){2}\put(-12,24){\vector(0,-1){22}}
\put(27,23){2}
\put(42,23){2} \put(36.5,23.2){$<$}

\put(3,17){$2\bar{3}$}
\put(19,17){2}
\put(34,17){2}

\put(-3,11){$\bar{3}$}\put(-5,12){\line(-1,0){4}}
\put(12,11){2}\put(-9,12){\vector(0,-1){10}}
\put(27,11){2}
\put(42,11){2}\put(36.5,11.2){$<$}

\put(3,5){22}
\put(19,5){2}
\put(34,5){2}

\put(-19,-1){33}\put(-12,-0.5){$\rightsquigarrow$}
\put(-4,-1){$\bar{2}\bar{2}$}
\put(11,-1){22}
\put(26,-1){22}
\put(41,-1){22}

\put(3,-7){$\bar{2}\bar{2}$}
\put(19,-7){2}
\put(34,-7){2}

\put(-4,-13){2}
\put(12,-13){$\bar{2}$}
\put(27,-13){2}
\put(42,-13){2}

\put(3,-19){$2\bar{2}$}
\put(19,-19){$\bar{2}$}
\put(34,-19){2}

\put(-4,-25){2}
\put(12,-25){$\bar{2}$}
\put(27,-25){$\bar{2}$}
\put(42,-25){2}

\put(3,-31){22}
\put(19,-31){$\bar{2}$}
\put(34,-31){$\bar{2}$}

\end{picture}
}
\end{equation*}
Thus we get $p^{(1)}=\Phi_2{\mathcal C}_2(p^{(2)}) = 
\boxed{22}\otimes 
\begin{picture}(12,1)
\put(0,0){$\boxed{\phantom{3}}$}
\put(3,0){$\bar{2}$}
\end{picture}
\otimes \begin{picture}(12,1)
\put(0,0){$\boxed{\phantom{3}}$}
\put(3,0){$\bar{2}$}
\end{picture}
\,$.
Assigning the mode and normal ordering, we find 
${\mathcal C}_1(p^{(1)}) = :
\boxed{22}_{\, 5}\otimes 
\begin{picture}(16,1)
\put(0,0){$\boxed{\phantom{3}}_{\, 1}$}
\put(3,0){$\bar{2}$}
\end{picture}
\otimes \begin{picture}(16,1)
\put(0,0){$\boxed{\phantom{3}}_{\, 4}$}
\put(3,0){$\bar{2}$}
\end{picture}
\,: = 
\boxed{2}_{\, 2} \otimes 
\begin{picture}(24,1)
\put(0,0){$\boxed{\phantom{33}}_{\, 4}$}
\put(3,0){$2\bar{2}$}
\end{picture} \otimes 
\begin{picture}(17,1)
\put(0,0){$\boxed{\phantom{3}}_{\, 4}$}
\put(3,0){$\bar{2}$}
\end{picture}
$.
Finally to find  
$\Phi_1{\mathcal C}_1(p^{(1)})$
we compute 
${\mathcal T}^2_1 \otimes \boxed{2} \otimes {\mathcal T}^2_1
\otimes 
\begin{picture}(18,1)
\put(0,0){$\boxed{\phantom{33}}$}
\put(3,0){$2\bar{2}$}
\end{picture} \otimes 
\begin{picture}(12,1)
\put(0,0){$\boxed{\phantom{3}}$}
\put(3,0){$\bar{2}$}
\end{picture}
\otimes \bigl(\;\boxed{1}^{\otimes 7}\,\bigr)$ by 
using the $A^{(2)}_5$ combinatorial $R$ on letters 
$1, 2, 3, \bar{3}, \bar{2}, \bar{1}$.
It is by given applying ${\mathcal T}_1^2\otimes$
to the output of the following diagram
($\,\boxed{1}^{\otimes 7}$ has been 
replaced with $\,\boxed{1}^{\otimes 5}$ below): 
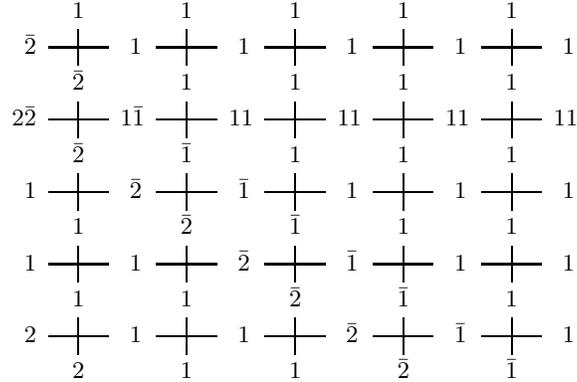
\begin{equation*}
\unitlength 0.8mm
{\small
\begin{picture}(110,65)(-10,-19)
\multiput(0,-12)(0,12){5}{
\multiput(0,0)(18,0){5}{\line(1,0){10}}
\multiput(5,-3)(18,0){5}{\line(0,1){6}}
}

\multiput(4,41)(18,0){5}{\put(0,0){1}}

\put(-4,35){$\bar{2}$}
\multiput(13.5,35)(18,0){5}{\put(0,0){1}}

\put(4,29){$\bar{2}$}
\multiput(22,29)(18,0){4}{\put(0,0){1}}

\put(-6,23){$2\bar{2}$}
\put(12,23){$1\bar{1}$}
\put(30,23){11}
\multiput(48,23)(18,0){3}{\put(0,0){11}}

\put(4,17){$\bar{2}$}\put(22,17){$\bar{1}$}\put(40,17){1}
\multiput(58,17)(18,0){2}{\put(0,0){1}}

\put(-4,11){$1$}\put(13.5,11){$\bar{2}$}\put(31.5,11){$\bar{1}$}
\put(49.5,11){1}
\multiput(67.5,11)(18,0){2}{\put(0,0){1}}

\put(4,5){1}\put(22,5){$\bar{2}$}\put(40,5){$\bar{1}$}
\put(58,5){1}\put(76,5){1}

\put(-4,-1){1}\put(13.5,-1){1}\put(31.5,-1){$\bar{2}$}
\put(49.5,-1){$\bar{1}$}\put(67.5,-1){1}
\put(85.5,-1){1}

\put(4,-7){1}\put(22,-7){1}
\put(40,-7){$\bar{2}$}
\put(58,-7){$\bar{1}$}\put(76,-7){1}

\put(-4,-13){2}\put(13.5,-13){1}\put(31.5,-13){1}
\put(49.5,-13){$\bar{2}$}\put(67.5,-13){$\bar{1}$}
\put(85.5,-13){1}

\put(4,-19){2}\put(22,-19){1}\put(40,-19){1}
\put(58,-19){$\bar{2}$}\put(76,-19){$\bar{1}$}

\end{picture}
}
\end{equation*}
Therefore we obtain 
$\Phi_1{\mathcal C}_1(p^{(1)})= 11211\bar{2}\bar{1}$ 
in agreement with (\ref{eq:pA(2)odd}).

\section{$B^{(1)}_n$ case}\label{sec:B(1)}

\subsection{Rigged configurations}

Consider the data of the form (\ref{eq:rc}),
where $\mu^{(a)}$ with $a<n$ and $2\mu^{(n)}$ are partitions.
Thus $\mu^{(n)}_i \in \Z/2$.
$J^{(a)} = (J^{(a)}_1,\ldots, J^{(a)}_{l_a})$ is taken from 
$(\Z_{\ge 0})^{l_a}$ for all $1 \le a \le n$.
Using $E^{(a)}_j$ in (\ref{eq:Eaj}),  we define the vacancy numbers
by (\ref{eq:paj}) with $1 \le a \le n-2, \, j \in \Z_{>0}$ and 
\begin{align*}
p^{(n-1)}_j &= E^{(n-2)}_j - 2E^{(n-1)}_j + 2E^{(n)}_j
\quad (j \in \Z_{>0}),\\
p^{(n)}_j &= 2E^{(n-1)}_j - 4E^{(n)}_j\quad (j \in \Z_{>0}/2).
\end{align*}
All the vacancy numbers are integers even though $E^{(n)}_j \in \Z/2$.
The data (\ref{eq:rc}) is called a $B^{(1)}_n$ rigged configuration
if (\ref{eq:cond}) is satisfied, 
where $j \in \Z_{>0}$ for $1 \le a \le n-1$ and 
$j \in \Z_{>0}/2$ for $a=n$.

The following is an example of $B^{(1)}_2$ rigged configuration.
\begin{equation}\label{eq:rcB}
{\small
\unitlength 10pt
\begin{picture}(20,6)(-2,-2)
\put(0,1){$(1^{7})$}
\put(0,3){$\mu^{(0)}$}

\multiput(1,0)(0,0){1}{
\put(4,1){0}
\put(4,-1){2}

\put(5,2){\line(1,0){2}}
\put(5,1){\line(1,0){2}}
\multiput(5,0)(0,-1){3}{\line(1,0){1}}
\put(5,2){\line(0,-1){4}}
\put(6,2){\line(0,-1){4}}
\put(7,2){\line(0,-1){1}}

\put(7.3,1.05){0}
\put(6.3,0.05){2}
\put(6.3,-0.95){2}
\put(6.3,-1.95){2}

\put(5,3){$\mu^{(1)}$}
}
\put(12,0){0}
\put(12,1){2}

\multiput(13,1)(0,1){2}{\line(1,0){1}}
\put(13,0){\line(1,0){0.5}}
\put(13,2){\line(0,-1){2}}
\put(13.5,2){\line(0,-1){2}}
\put(14,2){\line(0,-1){1}}

\put(13.8,0.05){0}
\put(14.3,1.05){1}
\put(13,3){$\mu^{(2)}$}
\end{picture}
}
\end{equation}
We have depicted $\mu^{(2)}$ by a Young diagram
consisting of $1\times \frac{1}{2}$ elementary blocks.

\subsection{Crystals}

\begin{align*}
B_l &=\{ (x_1,\ldots,x_n,x_0,\bar{x}_n,\ldots, \bar{x}_1)
\in (\Z_{\ge 0})^{2n+1} \mid 
x_0+\sum_{i=1}^n(x_i+\bar{x}_i) = l, x_0=0,1\},\\
B^{\ge a+1}_l &= \{ (x_1,\ldots,x_n,x_0,\bar{x}_n,\ldots, \bar{x}_1)
 \in B_l \mid 
x_i=\bar{x}_i=0 \, (1 \le i \le a)\}\quad (0 \le a \le n-1).
\end{align*}
\begin{equation*}
B_l = B^{\ge 1}_l \supset B^{\ge 2}_l \supset \cdots 
\supset B^{\ge n}_l 
\end{equation*}
as sets.
We regard $B^{\ge a+1}_l$ as $B^{(1)}_{n-a}$ crystal 
and represent the elements by 
length $l$ semistandard row tableaux over the alphabet 
$a+1 < \cdots <  n < 0 < \bar{n}<  \cdots < \overline{a+1}$.
We use the notation (\ref{eq:al}).
Then the highest element in $B^{\ge a+1}_l$ with respect to 
the classical part $B_{n-a}$ is $(a+1)^l$.

The $B^{(1)}_2$ rigged configuration (\ref{eq:rcB}) corresponds to 
\begin{equation}\label{eq:pB(1)}
\boxed{1} \otimes \boxed{1}\otimes \boxed{2}\otimes 
\boxed{2}\otimes \boxed{0}\otimes \boxed{0}\otimes \boxed{0} 
\,\in {\mathcal P}_+((1^7)).
\end{equation}

\subsection{\mathversion{bold} Maps 
${\mathcal C}_1, \ldots, {\mathcal C}_{n}$}

For $1 \le a \le n-1$, the map ${\mathcal C}_a$ is  
specified by (\ref{eq:Ca}) and (\ref{eq:di}) 
provided that $B_l$ is regarded as the 
$B^{(1)}_{n-a}$ crystal $B^{\ge a+1}_l$.

Define $\bar{p}^{(n)}$  and ${\mathcal C}_n({\bar p}^{(n)})$ 
in the same way as in Section \ref{subsec:cpcp}
by replacing $\mu^{(n)}_i$ with $2\mu^{(n)}_i$.

\subsection{\mathversion{bold} Maps 
$\Phi_1, \ldots, \Phi_{n}$}

For $1\le a \le n-1$, we define the map 
$\Phi_a$ by (\ref{eq:Phi}) and (\ref{eq:zdon}) 
by regarding $B_l = B^{\ge a+1}_l$ and 
$B'_l = B^{\ge a}_l$ as the $B^{(1)}_{n-a}$ and 
$B^{(1)}_{n-a+1}$ crystals.
Since the energy function can be negative we employ
the same regularization as $D^{(1)}_n$ explained in 
Section \ref{subsec:PhiD}.

To define $\Phi_n{\mathcal C}_n({\bar p}^{(n)})$, 
we put
$b_1[d_1]\otimes \cdots \otimes b_m[d_m] 
= {\mathcal C}_{n}({\bar p}^{(n)})$ and 
$(\lambda_1, \ldots, \lambda_k)=2\mu^{(n-1)}$ 
in  (\ref{eq:Phi}) 
regarding $B_l$ as the $A^{(1)}_0$ crystal with letter ${\bar n}$ only 
and $B'_l$ as the $A^{(1)}_1$ crystal with letters $n, {\bar n}$.
By using the resulting $c_1 \otimes \cdots \otimes c_k$ in (\ref{eq:Phi}), 
we construct $\Phi_n{\mathcal C}_n({\bar p}^{(n)})$ as 
\begin{equation}
\Phi_n{\mathcal C}_n({\bar p}^{(n)})
= {\tilde c}_1 \otimes \cdots \otimes {\tilde c}_k 
\in B^{\ge n}_{\lambda_1} \otimes \cdots \otimes B^{\ge n}_{\lambda_k}.
\end{equation}
Here for $B'_l$ with even $l$, $\tilde{c}_j$ is determined from
$c_j$ as
\begin{align}
\tilde{\phantom{c}}: \;
B'_l \qquad & \longrightarrow 
\quad \qquad \qquad B^{\ge n}_l \nonumber \\
(x_n, {\bar x}_n) \quad & \longmapsto 
(0,\ldots,0,z_{n}, z_0, \bar{z}_n, 0,\ldots,0)\nonumber\\
z_n = \left[\frac{{x}_n}{2}\right],  & 
{\bar z}_n = \left[\frac{{\bar x}_n}{2}\right], 
z_0 = l - z_n - {\bar z}_n,\label{eq:tilde}
\end{align}
where $[x]$ denotes the largest integer not exceeding $x$.

\subsection{Example}\label{subsec:exB}

Consider the $B^{(1)}_2$ rigged configuration (\ref{eq:rcB}).
Doubling the width of $\mu^{(2)}$,  we have 
\begin{equation*}
\bar{p}^{(2)} = 
\begin{picture}(18,1)
\put(0,0){$\boxed{\phantom{33}}$}
\put(3,0){$\bar{2}\bar{2}$}
\end{picture}
\otimes 
\begin{picture}(12,1)
\put(0,0){$\boxed{\phantom{2}}$}
\put(3,0){$\bar{2}$}
\end{picture},\quad 
{\mathcal C}_2(\bar{p}^{(2)}) =:
\begin{picture}(24,1)
\put(0,0){$\boxed{\phantom{33}}_{\, 3}$}
\put(3,0){$\bar{2}\bar{2}$}
\end{picture}
\otimes 
\begin{picture}(18,1)
\put(0,0){$\boxed{\phantom{2}}_{\, 2}$}
\put(3,0){$\bar{2}$}
\end{picture}:=
\begin{picture}(18,1)
\put(0,0){$\boxed{\phantom{2}}_{\, 1}$}
\put(3,0){$\bar{2}$}
\end{picture}\otimes 
\begin{picture}(24,1)
\put(0,0){$\boxed{\phantom{33}}_{\, 4}.$}
\put(3,0){$\bar{2}\bar{2}$}
\end{picture}
\end{equation*}
To find $\Phi_2{\mathcal C}_2(\bar{p}^{(2)})$ we first need 
${\mathcal T}_2^1\otimes 
\begin{picture}(12,1)
\put(0,0){$\boxed{\phantom{2}}$}
\put(3,0){$\bar{2}$}
\end{picture}
\otimes 
{\mathcal T}_2^3 \otimes 
\begin{picture}(18,1)
\put(0,0){$\boxed{\phantom{33}}$}
\put(3,0){$\bar{2}\bar{2}$}
\end{picture}
\otimes \bigl(\; \boxed{2222} \otimes 
\boxed{22} \otimes \boxed{22} \otimes \boxed{22} \,\bigr)$.
\begin{equation*}
\unitlength 0.8mm
{\small
\begin{picture}(110,77)(-20,-19)
\multiput(0,-12)(0,12){6}{
\multiput(0,0)(18,0){4}{\line(1,0){10}}
\multiput(5,-3)(18,0){4}{\line(0,1){6}}
}

\put(1,53){2222}\multiput(21,53)(18,0){3}{\put(0,0){22}}

\put(-6,47){$\bar{2}\bar{2}$}
\put(12,47){22}
\put(30,47){22}
\multiput(48,47)(18,0){2}{\put(0,0){22}}

\put(1,41){$22\bar{2}\bar{2}$}
\multiput(21,41)(18,0){3}{\put(0,0){22}}

\put(-5,35){2}\put(13,35){$\bar{2}$}
\multiput(31,35)(18,0){3}{\put(0,0){2}}

\put(1,29){$222\bar{2}$}\put(21,29){$2\bar{2}$}
\multiput(39,29)(18,0){2}{\put(0,0){22}}

\put(-5,23){2}\put(13,23){$\bar{2}$}
\multiput(31,23)(18,0){3}{\put(0,0){2}}

\put(1,17){2222}\put(21,17){$\bar{2}\bar{2}$}
\multiput(39,17)(18,0){2}{\put(0,0){22}}

\put(-5,11){2}\put(13,11){2}\put(31,11){$\bar{2}$}
\multiput(49,11)(18,0){2}{\put(0,0){2}}

\put(1,5){2222}\put(21,5){$2\bar{2}$}\put(39,5){$2\bar{2}$}
\put(57,5){22}

\put(-5,-1){$\bar{2}$}\put(13,-1){2}\put(31,-1){$\bar{2}$}
\multiput(49,-1)(18,0){2}{\put(0,0){2}}

\put(1,-7){$222\bar{2}$}\put(21,-7){22}\put(39,-7){$\bar{2}\bar{2}$}
\put(57,-7){22}

\put(-5,-13){2}\put(13,-13){$\bar{2}$}\put(31,-13){2}
\put(49,-13){$\bar{2}$}\put(67,-13){2}

\put(1,-19){2222}\multiput(21,-19)(18,0){3}{\put(0,0){$2\bar{2}$}}

\end{picture}
}
\end{equation*}
Under $\,\tilde{\phantom{2} }\,$, the elements $\boxed{2222}$ and 
$
\begin{picture}(18,1)
\put(0,0){$\boxed{\phantom{22}}$}
\put(3,0){$2\bar{2}$}
\end{picture}
$
turn into $\boxed{22}$ and $\boxed{0}$. 
Thus we get 
\begin{equation*}
p^{(1)}=\Phi_2{\mathcal C}_2(\bar{p}^{(2)}) = 
\boxed{22}\otimes \boxed{0}\otimes \boxed{0}\otimes \boxed{0},
\end{equation*}
which is an element of $B^{(1)}_1$ crystal with letters $2,0,\bar{2}$.
In view of  $\boxed{22}\otimes \boxed{0} \simeq 
\boxed{2} \otimes \boxed{20}$ and 
$H\bigl(\boxed{22}\otimes \boxed{0}\bigr) = 
H\bigl(\boxed{0}\otimes \boxed{0}\bigr) = 0$, 
we have ${\mathcal C}_1(p^{(1)}) = 
\boxed{22}_{\,2}\otimes 
\boxed{0}_{\,3}\otimes \boxed{0}_{\,3}\otimes \boxed{0}_{\,3}$.
$\Phi_1{\mathcal C}_1(p^{(1)})$ is derived from 
${\mathcal T}_1^2 \otimes \boxed{22} \otimes {\mathcal T}_1^1
\otimes \boxed{0}^{\otimes 3}\otimes \bigl(\; 
\boxed{1}^{\otimes 7}\, \bigr)$.
Since the combinatorial $R$ is the identity on $B_1 \otimes B_1$, 
the part ${\mathcal T}_1^1
\otimes \boxed{0}^{\otimes 3}$
just pushes $\boxed{1}$ to the right. 
$\Phi_1{\mathcal C}_1(p^{(1)})$ is obtained by 
putting $\boxed{1}^{\otimes 2}$ in front of the output of the 
diagram:
\begin{equation*}
\unitlength 0.8mm
{\small
\begin{picture}(110,17)(5,-19)
\multiput(0,-12)(0,12){1}{
\multiput(0,0)(18,0){7}{\line(1,0){10}}
\multiput(5,-3)(18,0){7}{\line(0,1){6}}
}

\put(4,-7){1}\put(22,-7){0}\put(40,-7){0}
\put(58,-7){0}\put(76,-7){1}\put(94,-7){1}\put(112,-7){1}

\put(-6,-13){22}\put(12,-13){12}\put(30,-13){10}
\put(48,-13){10}\put(66,-13){10}\put(84,-13){11}
\put(102,-13){11}\put(120,-13){11}

\put(4,-19){2}\put(22,-19){2}\put(40,-19){0}
\put(58,-19){0}\put(76,-19){0}\put(94,-19){1}\put(112,-19){1}
\end{picture}
}
\end{equation*}
Thus we get $\Phi_1{\mathcal C}_1(p^{(1)})=1122000$
in agreement with (\ref{eq:pB(1)}).

\section{$C^{(1)}_n$ case}\label{sec:C(1)}

\subsection{Rigged configurations}\label{subsec:rcC(1)}

Consider the data of the form (\ref{eq:rc}) and 
define $E^{(a)}_j$ as in (\ref{eq:Eaj}).
The vacancy number is specified by (\ref{eq:paj}) 
for $1 \le a \le n-1$ and 
\begin{equation}\label{eq:pajC}
p^{(n)}_j = E^{(n-1)}_j - E^{(n)}_j.
\end{equation}
The data (\ref{eq:rc}) is called a 
$C^{(1)}_n$ rigged configuration if 
$\mu^{(n)}_i \in 2\Z$ and (\ref{eq:cond}) is satisfied.
It is depicted as in $A^{(1)}_n$ case by the Young diagrams
with rigging.
The following is an example of $C^{(1)}_2$ rigged configuration.
\begin{equation}\label{eq:rcC}
{\small
\unitlength 10pt
\begin{picture}(20,6)(1,-2)
\put(0,1){$(1^{10})$}
\put(0,3){$\mu^{(0)}$}
\put(4,0){1}
\put(4,1){0}
\put(4,-1){6}
\multiput(5,0)(1,0){4}{\line(0,1){2}}
\put(5,0){\line(1,0){3}}
\put(5,-1){\line(1,0){1}}
\put(5,-1){\line(0,1){1}}
\put(6,-1){\line(0,1){1}}
\multiput(5,1)(0,1){2}{\line(1,0){4}}
\put(9,1){\line(0,1){1}}

\put(8.3,0.05){1}
\put(9.3,1.05){0}
\put(6.3,-0.95){2}
\put(6,3){$\mu^{(1)}$}
\put(12,0){1}
\put(12,1){2}
\multiput(13,0)(1,0){3}{\line(0,1){2}}
\put(13,0){\line(1,0){2}}
\multiput(13,1)(0,1){2}{\line(1,0){4}}
\put(16,1){\line(0,1){1}}
\put(17,1){\line(0,1){1}}
\put(15.3,0.05){0}
\put(17.3,1.05){2}
\put(13.7,3){$\mu^{(2)}$}
\end{picture}
}
\end{equation}

\subsection{\mathversion{bold} Crystals}

\begin{equation}
B_l =\{ (x_0,\ldots,x_n,\bar{x}_n,\ldots, \bar{x}_0)
\in (\Z_{\ge 0})^{2n+2} \mid 
\sum_{i=0}^n(x_i+\bar{x}_i) = l, x_0 = \bar{x}_0 \}.
\end{equation}
\begin{equation}
B^{\ge a}_l = \{ (x_0,\ldots,x_n,\bar{x}_n,\ldots, \bar{x}_0)
 \in B_l \mid 
x_i=\bar{x}_i=0 \, (0 \le i \le a-1), 
x_a=\bar{x}_a \}
\end{equation}
for $0 \le a \le n$.
Then 
\begin{equation}\label{eq:uniqueC}
B_l = B^{\ge 0}_l \supset B^{\ge 1}_l \supset \cdots 
\supset B^{\ge n}_l =
\begin{cases} 
\emptyset & l : {\rm odd}\\
\{(0,\ldots, 0, \frac{l}{2}, \frac{l}{2},0, \ldots,0)\} & l: {\rm even}
\end{cases}
\end{equation}
as sets.
We regard $B^{\ge a}_l$ as $C^{(1)}_{n-a}$ crystal 
and represent the elements by 
length $l$ semistandard row tableaux over the alphabet 
$a < \cdots <  n < \bar{n}<  \cdots < \bar{a}$.
For example for $n=3, a=2$, 
the highest elements in $B^{\ge 2}_5$ with respect to 
$C_{n-a}$ 
are
\begin{equation*}
\boxed{33333}, \;
\boxed{2333\bar{2}},\;
\boxed{223\bar{2}\bar{2}}.
\end{equation*}
The rigged configuration (\ref{eq:rcC}) corresponds to
\begin{equation}\label{eq:Cpath}
\unitlength 10pt
\begin{picture}(21,1.1)(1,0.1)
\multiput(-0.3,0.1)(2.3,0){10}{$\boxed{\phantom{3}}$}
\multiput(1.1,0.1)(2.3,0){9}{$\otimes$}
\put(0,0){1}
\put(2.3,0){$1$}
\put(4.6,0){$1$}
\put(6.9,0){$2$}
\put(9.2,0){$\bar{2}$}
\put(11.5,0){$2$}
\put(13.8,0){$2$}
\put(16.1,0){$1$}
\put(18.4,0){${\bar 1}$}
\put(20.7,0){${\bar 1}$}
\end{picture}\in {\mathcal P}_+((1^{10})).
\end{equation}

\subsection{\mathversion{bold} Maps ${\mathcal C}_1, 
\ldots, {\mathcal C}_n$}

Let $B_l = B^{\ge a}_l$ be the $C^{(1)}_{n-a}$ crystal.
Then the map ${\mathcal C}_a$ with $a<n$ 
is defined in the same way as in 
Section \ref{subsec:clclA(1)}.
Since the energy function is nonnegative,
one has $0 \le d_1 \le \cdots \le d_m$ for 
normal ordered forms
$b_1[d_1]\ot \cdots \ot b_m[d_m]$.
Define the $C^{(1)}_0$ crystal element 
\begin{equation}\label{eq:pnC}
\bar{p}^{(n)} = (n\bar{n})^{\mu^{(n)}_1/2}\ot \cdots \ot 
(n\bar{n})^{\mu^{(n)}_{l_n}/2},
\end{equation}
where $(n\bar{n})^{l/2}$ denotes the unique
element in $B^{\ge n}_l$ with even $l$ in (\ref{eq:uniqueC}).
(Recall that $\mu^{(n)}_i$ is even for $C^{(1)}_n$ rigged configuration.)
Then ${\mathcal C}_n(\bar{p}^{(n)})$
is obtained by assigning the mode in (\ref{eq:Ca}) as
$d_i = J_i + \mu_i/2$ instead of (\ref{eq:di}).
For normal ordering we need the 
$C^{(1)}_0$ combinatorial $R$.
It is formally specified by
saying that $H=0$ on $B_l \otimes B_m$ and 
$b[d]\otimes b'[d'] \simeq b'[d']\otimes b[d]$.

\subsection{\mathversion{bold} Maps $\Phi_1, \ldots, \Phi_n$}

The map $\Phi_a$ is defined in the same way as in 
Section \ref{subsec:Phi} if 
$B_l$ and $B'_l$ there are replaced by 
$B_l = B^{\ge a}_l$ and $B'_l = B^{\ge a-1}_l$, respectively.

\subsection{Example}

We start from the $C^{(1)}_2$ rigged configuration  (\ref{eq:rcC}).
Then 
\begin{equation*}
\bar{p}^{(2)} = 
\begin{picture}(28,1.2)
\put(0,0){$\boxed{\phantom{2222}}$}
\put(3,0){$22\bar{2}\bar{2}$}
\end{picture}
\otimes 
\begin{picture}(18,1.2)
\put(0,0){$\boxed{\phantom{22}},$}
\put(3,0){$2\bar{2}$}
\end{picture}
\qquad
{\mathcal C}_2(\bar{p}^{(2)}) =
\begin{picture}(18,1.2)
\put(0,0){$\boxed{\phantom{22}}_{\, 1}$}
\put(3,0){$2\bar{2}$}
\end{picture} 
\;\;\, \otimes 
\begin{picture}(28,1.2)
\put(0,0){$\boxed{\phantom{2222}}_{\, 4}.$}
\put(3,0){$22\bar{2}\bar{2}$}
\end{picture}
\end{equation*}
In order to find $\Phi_2{\mathcal C}_2(\bar{p}^{(2)})$ we are to compute 
\begin{equation*}
{\mathcal T}_2^1 \otimes 
\begin{picture}(18,1.2)
\put(0,0){$\boxed{\phantom{22}}$}
\put(3,0){$2\bar{2}$}
\end{picture}
\otimes 
{\mathcal T}_2^3
\otimes 
\begin{picture}(28,1.2)
\put(0,0){$\boxed{\phantom{2222}}$}
\put(3,0){$22\bar{2}\bar{2}$}
\end{picture}
\otimes \bigl(
\boxed{2222} \otimes \boxed{222} \otimes \boxed{2}
\bigr)
\end{equation*}
for $C^{(1)}_1$ crystal with letters $1,2,\bar{2},\bar{1}$.
{}From the diagram
\begin{equation*}
\unitlength 1mm
\begin{picture}(30,77)(5,-32)
\multiput(0,-24)(0,12){6}{
\multiput(0,0)(18,0){3}{\line(1,0){10}}
\multiput(5,-3)(18,0){3}{\line(0,1){6}}
}

\put(2,41){2222}
\put(20.5,41){222}
\put(40,41){2}

\put(-8,35){$22\bar{2}\bar{2}$}
\put(11,35){$2222$}
\put(29,35){$2222$}
\put(47,35){$2222$}

\put(2,29){$22\bar{2}\bar{2}$}
\put(20.5,29){222}
\put(40,29){2}

\put(-5,23){2}
\put(13,23){2}
\put(31,23){2}
\put(49,23){2}

\put(2,17){$12\bar{2}\bar{1}$}
\put(20.5,17){222}
\put(40,17){2}

\put(-5,11){2}
\put(13,11){2}
\put(31,11){2}
\put(49,11){2}

\put(2,5){$11\bar{1}\bar{1}$}
\put(20.5,5){222}
\put(40,5){2}

\put(-5,-1){2}
\put(13,-1){$\bar{2}$}
\put(31,-1){2}
\put(49,-1){2}

\put(2,-7){$122\bar{1}$}
\put(20.5,-7){$22\bar{2}$}
\put(40,-7){2}

\put(-5,-13){$2\bar{2}$}
\put(13,-13){$1\bar{1}$}
\put(31,-13){22}
\put(49,-13){22}

\put(2,-19){$222\bar{2}$}
\put(20.5,-19){$1\bar{2}\bar{1}$}
\put(40,-19){2}

\put(-5,-25){2}
\put(13,-25){2}
\put(31,-25){$\bar{2}$}
\put(49,-25){2}

\put(2,-31){$122\bar{1}$}
\put(20.5,-31){$12\bar{1}$}
\put(40,-31){$\bar{2}$}

\end{picture}
\end{equation*}
we get $p^{(1)} = \Phi_2{\mathcal C}_2(\bar{p}^{(2)})=
\begin{picture}(28,1.2)
\put(0,0){$\boxed{\phantom{2222}}$}
\put(3,0){$122\bar{1}$}
\end{picture}
\otimes 
\begin{picture}(23,1.2)
\put(0,0){$\boxed{\phantom{122}}$}
\put(3,0){$12\bar{1}$}
\end{picture}
\otimes 
\begin{picture}(13,1.2)
\put(0,0){$\boxed{\phantom{2}}$}
\put(3,0){$\bar{2}$}
\end{picture}
$.
Assigning the mode to $p^{(1)}$ we have 
\vspace{0.2cm}
\\
${\mathcal C}_1(p^{(1)}) = 
\begin{picture}(40,1.2)
\put(0,0){$:\, \boxed{\phantom{2222}}_{\, 3}$}
\put(10,0){$122\bar{1}$}
\end{picture}
\otimes 
\begin{picture}(29,1.2)
\put(0,0){$\boxed{\phantom{122}}_{\, 4}$}
\put(3,0){$12\bar{1}$}
\end{picture}
\otimes 
\begin{picture}(13,1.2)
\put(0,0){$\boxed{\phantom{2}}_{\, 3} :$}
\put(4,0){$\bar{2}$}
\end{picture}
\quad$, where the normal ordering is to be done along the 
$C^{(1)}_1$ crystal with letters $1,2,\bar{2},\bar{1}$.
There are four normal ordered forms:
\begin{align*}
&
\begin{picture}(18,1.2)
\put(0,0){$\boxed{\phantom{2}}_{\, 3}$}
\put(4,0){$2$}
\end{picture}
\otimes
\begin{picture}(29,1.2)
\put(0,0){$\boxed{\phantom{122}}_{\, 3}$}
\put(3,0){$22\bar{2}$}
\end{picture}
\otimes 
\begin{picture}(35,1.2)
\put(0,0){$\boxed{\phantom{2222}}_{\, 4}$}
\put(3,0){$11\bar{1}\bar{1}$}
\end{picture},\quad 
\begin{picture}(20,1.2)
\put(0,0){$\boxed{\phantom{2}}_{\, 3}$}
\put(4,0){$2$}
\end{picture}
\otimes 
\begin{picture}(35,1.2)
\put(0,0){$\boxed{\phantom{2222}}_{\, 3}$}
\put(3,0){$222\bar{2}$}
\end{picture}
\otimes
\begin{picture}(29,1.2)
\put(0,0){$\boxed{\phantom{122}}_{\, 4}$}
\put(3,0){$1\bar{2}\bar{1}$}
\end{picture}\\
& \\&
\begin{picture}(29,1.2)
\put(0,0){$\boxed{\phantom{122}}_{\, 3}$}
\put(3,0){$12\bar{1}$}
\end{picture}
\otimes
\begin{picture}(18,1.2)
\put(0,0){$\boxed{\phantom{2}}_{\, 3}$}
\put(3,0){$2$}
\end{picture}
\otimes
\begin{picture}(35,1.2)
\put(0,0){$\boxed{\phantom{2222}}_{\, 4}$}
\put(3,0){$11\bar{1}\bar{1}$}
\end{picture}
,\quad 
\begin{picture}(35,1.2)
\put(0,0){$\boxed{\phantom{2222}}_{\, 3}$}
\put(3,0){$122\bar{1}$}
\end{picture}
\otimes
\begin{picture}(18,1.2)
\put(0,0){$\boxed{\phantom{2}}_{\, 3}$}
\put(4,0){$2$}
\end{picture}
\otimes
\begin{picture}(29,1.2)
\put(0,0){$\boxed{\phantom{122}}_{\, 4}$}
\put(3,0){$1\bar{2}\bar{1}$}
\end{picture}
.
\end{align*}
We pick the first one as ${\mathcal C}_1(p^{(1)})$. 
Then $\Phi_1{\mathcal C}_1(p^{(1)})$ 
is calculated {}from 
\begin{equation*}
{\mathcal T}_1^3 \otimes 
\begin{picture}(13,1.2)
\put(0,0){$\boxed{2}$}
\end{picture}
\otimes
\begin{picture}(23,1.2)
\put(0,0){$\boxed{\phantom{222}}$}
\put(3,0){$22\bar{2}$}
\end{picture}
\otimes 
{\mathcal T}_1^1
\otimes 
\begin{picture}(28,1.2)
\put(0,0){$\boxed{\phantom{2222}}$}
\put(3,0){$11\bar{1}\bar{1}$}
\end{picture}
\otimes \bigl(\;\;
\boxed{1}^{\otimes 10} \;\bigr)
\end{equation*}
for $C^{(1)}_2$ crystal with letters $0,1,2,\bar{2},\bar{1},\bar{0}$.
It is given by applying ${\mathcal T}_1^3\otimes$
to the output of the following diagram
($\,\boxed{1}^{\otimes 10}$ has been 
replaced with $\,\boxed{1}^{\otimes 7}$ below): 
\begin{equation*}
\unitlength 0.8mm
{\small
\begin{picture}(110,53)(10,-19)
\multiput(0,-12)(0,12){4}{
\multiput(0,0)(18,0){7}{\line(1,0){10}}
\multiput(5,-3)(18,0){7}{\line(0,1){6}}
}

\multiput(4,29)(18,0){7}{\put(0,0){1}}

\put(-8.2,23){$11\bar{1}\bar{1}$}
\put(10.2,23){$111\bar{1}$}
\multiput(28.2,23)(18,0){6}{\put(0,0){1111}}

\put(4,17){$\bar{1}$}\put(22,17){$\bar{1}$}
\multiput(40,17)(18,0){5}{\put(0,0){1}}

\put(-4,11){$1$}\put(13.5,11){$\bar{1}$}\put(31.5,11){$\bar{1}$}
\multiput(49.5,11)(18,0){5}{\put(0,0){1}}

\put(4,5){$1$}\put(22,5){$\bar{1}$}\put(40,5){$\bar{1}$}
\multiput(58,5)(18,0){4}{\put(0,0){1}}

\put(-7,-1){$22\bar{2}$}\put(11,-1){122}\put(29,-1){$02\bar{0}$}
\put(47,-1){$0\bar{1}\bar{0}$}\put(65,-1){$1\bar{1}\bar{1}$}
\put(83,-1){$11\bar{1}$}
\multiput(101.5,-1)(18,0){2}{\put(0,0){111}}

\put(4,-7){$\bar{2}$}\put(22,-7){2}\put(40,-7){2}\put(58,-7){1}
\put(76,-7){$\bar{1}$}\put(94,-7){$\bar{1}$} \put(112,-7){1}

\put(-4,-13){$2$}\put(13.5,-13){$\bar{2}$}\put(31.5,-13){$2$}
\put(49.5,-13){2}\put(67.5,-13){1}
\put(85.5,-13){$\bar{1}$}\put(103.5,-13){$\bar{1}$}\put(121.5,-13){1}

\put(4,-19){2}\put(22,-19){$\bar{2}$}\put(40,-19){2}
\put(58,-19){2}\put(76,-19){1}\put(94,-19){$\bar{1}$}
\put(112,-19){$\bar{1}$}

\end{picture}
}
\end{equation*}
Thus we find 
$\Phi_1{\mathcal C}_1\Phi_2{\mathcal C}_2(\bar{p}^{(2)})
=1112\bar{2}221\bar{1}\bar{1}$ in agreement with 
(\ref{eq:Cpath}).

\section{$A^{(2)}_{2n}$ and $D^{(2)}_{n+1}$ cases}
\label{sec:III}

In this section we treat 
the algebras $\geh_n=A^{(2)}_{2n}$ and $D^{(2)}_{n+1}$ 
simultaneously.
They are reduced to the $C^{(1)}_n$ case.

\subsection{Rigged configurations}\label{subsec:rcIII}

Consider the data of the form (\ref{eq:rc}) and 
define $E^{(a)}_j$ as in (\ref{eq:Eaj}).
The vacancy number is specified by (\ref{eq:paj}) 
for $1 \le a \le n-1$ and 
\begin{align}
p^{(n)}_j &= E^{(n-1)}_j - E^{(n)}_j \quad\quad 
\hbox{ for } \; A^{(2)}_{2n},\label{eq:pajA(2)even}\\
p^{(n)}_j &= 2E^{(n-1)}_j - 2E^{(n)}_j \quad
\hbox{ for } \; D^{(2)}_{n+1}.\label{eq:pajD(2)}
\end{align}
The data (\ref{eq:rc}) is called a 
$\geh_n$ rigged configuration if 
(\ref{eq:cond}) is satisfied.
It is depicted as in $A^{(1)}_n$ case by the Young diagrams
with rigging.
For example 
\begin{equation}\label{eq:rcA(2)even}
{\small
\unitlength 10pt
\begin{picture}(20,6)(1,-2)
\put(0,1){$(1^{8})$}
\put(0,3){$\mu^{(0)}$}
\put(4,0){3}
\put(4,1){0}

\multiput(5,0)(1,0){3}{\line(0,1){2}}
\put(5,2){\line(1,0){4}}
\put(5,1){\line(1,0){4}}
\put(5,0){\line(1,0){2}}
\put(8,1){\line(0,1){1}}
\put(9,1){\line(0,1){1}}

\put(7.3,0.05){2}
\put(9.3,1.1){0}

\put(6,3){$\mu^{(1)}$}
\put(12,0){0}
\put(12,1){1}
\multiput(13,0)(1,0){2}{\line(0,1){2}}
\put(13,0){\line(1,0){1}}
\multiput(13,1)(0,1){2}{\line(1,0){3}}
\put(15,1){\line(0,1){1}}
\put(16,1){\line(0,1){1}}

\put(14.3,0.05){0}
\put(16.3,1.05){0}
\put(13.7,3){$\mu^{(2)}$}
\end{picture}
}
\end{equation}
is an $A^{(2)}_4$ rigged configuration, and 
\begin{equation}\label{eq:rcD(2)}
{\small 
\unitlength 10pt
\begin{picture}(20,6)(1,-2)
\put(0,1.1){$(1^{10})$}
\put(0,3.1){$\mu^{(0)}$}
\put(4,0){1}
\put(4,-1.9){5}
\put(5,2){\line(0,-1){4}}\put(6,2){\line(0,-1){4}}\put(7,2){\line(0,-1){3}}
\multiput(5,2)(0,-1){4}{\put(0,0){\line(1,0){2}}}
\put(5,-2){\line(1,0){1}}

\put(7.3,1.1){0}
\put(7.3,0.1){1}
\put(7.3,-0.9){1}
\put(6.3,-1.9){2}

\put(6,3){$\mu^{(1)}$}
\put(12,1){2}
\put(12,0){4}
\put(12,-1){2}

\put(13,2){\line(0,-1){3}}\put(14,2){\line(0,-1){3}}
\put(15,2){\line(0,-1){2}}\put(16,2){\line(0,-1){1}}
\put(13,2){\line(1,0){3}}\put(13,1){\line(1,0){3}}
\put(13,0){\line(1,0){2}}\put(13,-1){\line(1,0){1}}

\put(16.3,1.1){0}
\put(15.3,0.1){3}
\put(14.3,-0.9){1}

\put(13.7,3){$\mu^{(2)}$}
\end{picture}
}
\end{equation}
is a $D^{(2)}_{3}$ rigged configuration.
Let ${\rm RC}(\lambda)$ be the set of 
$\geh_n$ rigged configurations 
(\ref{eq:rc}) with $\mu^{(0)} = \lambda$.
One has the embedding:
{\small 
\begin{equation}\label{eq:embIII}
\begin{split}
\iota: {\rm RC}(\lambda) \hbox{ for } \geh_n \qquad\quad
&\longrightarrow \qquad \quad 
{\rm RC}(2\lambda) \hbox{ for } C^{(1)}_{n}\\
(\lambda, (\mu^{(1)},J^{(1)}), \ldots, (\mu^{(n)},J^{(n)}))
& \mapsto
\begin{cases}
(2\lambda, (2\mu^{(1)},2J^{(1)}), \ldots, (2\mu^{(n)}, 2J^{(n)}))
\hbox{ for }A^{(2)}_{2n},\\
(2\lambda, (2\mu^{(1)},2J^{(1)}), \ldots, (2\mu^{(n)},  J^{(n)}))
\hbox{ for }D^{(2)}_{n+1},
\end{cases}
\end{split}
\end{equation}
}
where for an array $\lambda=(\lambda_i)$, 
$2\lambda$ means $(2\lambda_i)$.
In the $D^{(2)}_{n+1}$ case, all the rigging except the $n$-th one 
$J^{(n)}$ are doubled.

\subsection{\mathversion{bold} Crystals}

For $A^{(2)}_{2n}$,
\begin{equation}\label{eq:BA(2)even}
B_l = \{(x_1,\ldots,x_n, x_\emptyset, \bar{x}_n,\ldots, \bar{x}_1)
\in (\Z_{\ge 0})^{2n+1} \mid 
x_\emptyset + \sum_{i=1}^n(x_i+\bar{x}_i) = l\}.
\end{equation}
For $D^{(2)}_{n+1}$,
\begin{equation}\label{eq:BD(2)}
B_l=\{(x_1, \ldots,x_n, x_0, x_\emptyset, \bar{x}_n,\ldots, \bar{x}_1)
\in (\Z_{\ge 0})^{2n+2} \mid 
x_0+x_\emptyset + \sum_{i=1}^n(x_i+\bar{x}_i) = l,
x_0 =0,1\}.
\end{equation}
There is an embedding
$(B_l \hbox{ for } \geh_n) \rightarrow 
(B_{2l} \hbox{ for } C^{(1)}_n)$
as sets, which will also be denoted by $\iota$:
\begin{equation}\label{eq:iota}
\iota(x)=\begin{cases}
(x_\emptyset, 2x_1,\ldots,2x_n,2\bar{x}_n,\ldots, 2\bar{x}_1,
x_\emptyset) &\hbox{for } A^{(2)}_{2n},\\
(x_\emptyset, 2x_1,\ldots,2x_{n-1},2x_n+x_0, 2\bar{x}_n+x_0,
2\bar{x}_{n-1},\ldots, 2\bar{x}_1, x_\emptyset) 
&\hbox{for } D^{(2)}_{n+1}.
\end{cases}
\end{equation}
Here $x =(x_1,\ldots, \bar{x}_1) 
\in B_l$ is specified as in (\ref{eq:BA(2)even}) and 
(\ref{eq:BD(2)}).
When applying $\iota^{-1}$, the tableau letters are 
halved (in case the image exists). 
In particular, a pair of 0 and $\bar{0}$ turns into 
$\emptyset$. As for $n$ and $\bar{n}$, the change is 
described by (\ref{eq:tilde}).
We extend the map $\iota$ to 
the tensor product of $\geh_n$ crystals by
$\iota(p_1 \otimes \cdots \otimes p_k) = 
\iota(p_1)\otimes \cdots \otimes \iota(p_k)$.

The $A^{(2)}_4$ rigged configuration (\ref{eq:rcA(2)even}) 
corresponds to 
\begin{equation}\label{eq:pA(2)even}
\boxed{1}\otimes \boxed{1}\otimes 
\begin{picture}(12,1)
\put(0,0){$\boxed{\phantom{3}}$}
\put(3,0){$\emptyset$}
\end{picture}
\otimes 
\begin{picture}(12,1)
\put(0,0){$\boxed{\phantom{3}}$}
\put(3,0){$\bar{1}$}
\end{picture}
\otimes \boxed{1}\otimes 
\begin{picture}(12,1)
\put(0,0){$\boxed{\phantom{3}}$}
\put(3,0){$\emptyset$}
\end{picture}\otimes \boxed{2}\otimes \boxed{2}\,
\in {\mathcal P}_+((1^8)).
\end{equation}
The $D^{(2)}_3$ rigged configuration (\ref{eq:rcD(2)}) 
corresponds to 
\begin{equation}\label{eq:pD(2)}
\boxed{1}\otimes 
\begin{picture}(12,1)
\put(0,0){$\boxed{\phantom{3}}$}
\put(3,0){$\emptyset$}
\end{picture}
\otimes \boxed{2} \otimes \boxed{1}\otimes \boxed{0} 
\otimes \boxed{1} \otimes 
\begin{picture}(12,1)
\put(0,0){$\boxed{\phantom{3}}$}
\put(3,0){$\bar{2}$}
\end{picture}
\otimes \boxed{2}\otimes 
\begin{picture}(12,1)
\put(0,0){$\boxed{\phantom{3}}$}
\put(3,0){$\emptyset$}
\end{picture}\otimes \boxed{0}\,
\in {\mathcal P}_+((1^{10})).
\end{equation}

\subsection{Example}\label{subsec:exIII}

Consider the $A^{(2)}_4$ rigged configuration 
(\ref{eq:rcA(2)even}).
According to (\ref{eq:embIII}) we double everything 
horizontally to get
\begin{equation}\label{eq:rrcA(2)even}
{\small 
\unitlength 9pt
\begin{picture}(20,6)(1,-1)
\put(0,2){\line(1,0){2}}\put(0,1){\line(1,0){2}}
\put(0,2){\line(0,-1){1}}\put(1,2){\line(0,-1){1}}\put(2,2){\line(0,-1){1}}
\put(2.2,2.0){8}

\put(0,3){$\mu^{(0)}$}
\put(4,1){0}
\put(4,0){6}
\put(5,2){\line(0,-1){2}}\put(6,2){\line(0,-1){2}}
\put(7,2){\line(0,-1){2}}\put(8,2){\line(0,-1){2}}
\put(9,2){\line(0,-1){2}}
\multiput(10,2)(1,0){4}{\put(0,0){\line(0,-1){1}}}

\multiput(5,2)(0,-1){2}{\put(0,0){\line(1,0){8}}}
\put(5,0){\line(1,0){4}}

\put(13.3,1.1){0}
\put(9.3,0.1){4}

\put(8,3){$\mu^{(1)}$}

\put(15,1.1){2}
\put(15,0.1){0}

\put(16,2){\line(0,-1){2}}\put(17,2){\line(0,-1){2}}
\put(18,2){\line(0,-1){2}}
\multiput(19,2)(1,0){4}{\put(0,0){\line(0,-1){1}}}
\put(16,2){\line(1,0){6}}
\put(16,1){\line(1,0){6}}
\put(16,0){\line(1,0){2}}

\put(22.3,1.1){0}
\put(18.3,0.1){0}

\put(17.5,3){$\mu^{(2)}$}
\end{picture}
}
\end{equation}
We regard this as a $C^{(1)}_2$ rigged configuration.
The vacancy numbers are 
determined by (\ref{eq:pajC}) which is the same as 
(\ref{eq:pajA(2)even}).
As the result they are all doubled as well. 
{}From (\ref{eq:pnC}) we set
\begin{equation*}
\bar{p}^{(2)} = 
\begin{picture}(38,1)
\put(0,0){$\boxed{\phantom{222222}}$}
\put(3,0){$222\bar{2}\bar{2}\bar{2}$}
\end{picture}
\otimes 
\begin{picture}(18,1)
\put(0,0){$\boxed{\phantom{22}}$}
\put(3,0){$2\bar{2}$}
\end{picture},\quad 
{\mathcal C}_2(\bar{p}^{(2)}) = 
\begin{picture}(22,1)
\put(0,0){$\boxed{\phantom{22}}_{\, 1}$}
\put(3,0){$2\bar{2}$}
\end{picture}
\otimes 
\begin{picture}(48,1)
\put(0,0){$\boxed{\phantom{222222}}_{\, 3}.$}
\put(3,0){$222\bar{2}\bar{2}\bar{2}$}
\end{picture}
\end{equation*}
Then $p^{(1)}:=\Phi_2{\mathcal C}_2(\bar{p}^{(2)})$ is determined as
\begin{equation*}
p^{(1)}= 
\begin{picture}(48,1)
\put(0,0){$\boxed{\phantom{22222222}}$}
\put(3,0){$11122\bar{1}\bar{1}\bar{1}$}
\end{picture}
\otimes 
\begin{picture}(28,1)
\put(0,0){$\boxed{\phantom{2222}}$}
\put(3,0){$122\bar{1}$}
\end{picture},\quad 
{\mathcal C}_1(p^{(1)}) = \begin{picture}(52,1)
\put(0,0){$\boxed{\phantom{22222222}}_{\, 5}$}
\put(3,0){$11122\bar{1}\bar{1}\bar{1}$}
\end{picture}
\otimes 
\begin{picture}(32,1)
\put(0,0){$\boxed{\phantom{2222}}_{\, 7}.$}
\put(3,0){$122\bar{1}$}
\end{picture}
\end{equation*}
$p=\Phi_1{\mathcal C}_1(p^{(1)})$ reads
\begin{equation*}
p=\boxed{11}\otimes \boxed{11} \otimes 
\begin{picture}(18,1)
\put(0,0){$\boxed{\phantom{22}}$}
\put(3,0){$0\bar{0}$}
\end{picture}\otimes
\begin{picture}(18,1)
\put(0,0){$\boxed{\phantom{22}}$}
\put(3,0){$\bar{1}\bar{1}$}
\end{picture}\otimes 
\boxed{11} \otimes 
\begin{picture}(18,1)
\put(0,0){$\boxed{\phantom{22}}$}
\put(3,0){$0\bar{0}$}
\end{picture}\otimes \boxed{22} \otimes \boxed{22}\,.
\end{equation*}
Finally applying $\iota^{-1}$ in (\ref{eq:iota}), we 
find $\iota^{-1}(p) = 
1\, 1\, \emptyset\, \bar{1}\, 1\, \emptyset\,  2\, 2, $
reproducing (\ref{eq:pA(2)even}).

Next we consider the $D^{(2)}_3$ rigged configuration (\ref{eq:rcD(2)}).
According to (\ref{eq:embIII}) we double everything horizontally 
except the rigging attached to $\mu^{(2)}$:
\begin{equation}\label{eq:rrcD(2)}
{\small 
\unitlength 9pt
\begin{picture}(20,6)(1,-2)
\put(0,2){\line(1,0){2}}\put(0,1){\line(1,0){2}}
\put(0,2){\line(0,-1){1}}\put(1,2){\line(0,-1){1}}\put(2,2){\line(0,-1){1}}
\put(2.2,2.0){10}

\put(0,3.1){$\mu^{(0)}$}
\put(4,0){2}
\put(3.5,-1.9){10}
\put(5,2){\line(0,-1){4}}\put(6,2){\line(0,-1){4}}
\put(7,2){\line(0,-1){4}}\put(8,2){\line(0,-1){3}}\put(9,2){\line(0,-1){3}}
\multiput(5,2)(0,-1){4}{\put(0,0){\line(1,0){4}}}
\put(5,-2){\line(1,0){2}}

\put(9.3,1.1){0}
\put(9.3,0.1){2}
\put(9.3,-0.9){2}
\put(7.3,-1.9){4}

\put(6.5,3){$\mu^{(1)}$}
\put(12,1){2}
\put(12,0){4}
\put(12,-1){2}

\put(13,2){\line(0,-1){3}}\put(14,2){\line(0,-1){3}}
\put(15,2){\line(0,-1){3}}\put(16,2){\line(0,-1){2}}
\put(17,2){\line(0,-1){2}}\put(18,2){\line(0,-1){1}}\put(19,2){\line(0,-1){1}}
\put(13,2){\line(1,0){6}}\put(13,1){\line(1,0){6}}
\put(13,0){\line(1,0){4}}\put(13,-1){\line(1,0){2}}

\put(19.3,1.1){0}
\put(17.3,0){3}
\put(15.3,-1){1}

\put(14.7,3){$\mu^{(2)}$}
\end{picture}
}
\end{equation}
We regard this as a $C^{(1)}_2$ rigged configuration.
The vacancy numbers 
determined by (\ref{eq:pajC}) instead of (\ref{eq:pajD(2)}) 
have been doubled except $p^{(2)}_j$'s 
which remain unchanged. 
{}From (\ref{eq:pnC}) we set
\begin{equation*}
\bar{p}^{(2)} = 
\begin{picture}(38,1)
\put(0,0){$\boxed{\phantom{222222}}$}
\put(3,0){$222\bar{2}\bar{2}\bar{2}$}
\end{picture}
\otimes 
\begin{picture}(28,1)
\put(0,0){$\boxed{\phantom{2222}}$}
\put(3,0){$22\bar{2}\bar{2}$}
\end{picture}
\otimes
\begin{picture}(18,1)
\put(0,0){$\boxed{\phantom{22}}$}
\put(3,0){$2\bar{2}$}
\end{picture},
\quad 
{\mathcal C}_2(\bar{p}^{(2)}) = 
\begin{picture}(22,1)
\put(0,0){$\boxed{\phantom{22}}_{\, 2}$}
\put(3,0){$2\bar{2}$}
\end{picture}
\otimes 
\begin{picture}(44,1)
\put(0,0){$\boxed{\phantom{222222}}_{\, 3}$}
\put(3,0){$222\bar{2}\bar{2}\bar{2}$}
\end{picture}
\otimes 
\begin{picture}(34,1)
\put(0,0){$\boxed{\phantom{2222}}_{\, 5}\,.$}
\put(3,0){$22\bar{2}\bar{2}$}
\end{picture}
\end{equation*}
Then $p^{(1)}=\Phi_2{\mathcal C}_2(\bar{p}^{(2)})$ is determined as
\begin{equation*}
p^{(1)} = \begin{picture}(28,1)
\put(0,0){$\boxed{\phantom{2222}}$}
\put(3,0){$12\bar{2}\bar{1}$}
\end{picture} \otimes 
\begin{picture}(28,1)
\put(0,0){$\boxed{\phantom{2222}}$}
\put(3,0){$22\bar{2}\bar{2}$}
\end{picture} \otimes 
\begin{picture}(28,1)
\put(0,0){$\boxed{\phantom{2222}}$}
\put(3,0){$12\bar{2}\bar{1}$}
\end{picture} \otimes 
\begin{picture}(18,1)
\put(0,0){$\boxed{\phantom{22}}$}
\put(3,0){$2\bar{2}$}
\end{picture}.
\end{equation*}
There are two normal ordered forms for 
${\mathcal C}_1(p^{(1)})$:
\begin{equation*}
\begin{picture}(30,1)
\put(0,0){$\boxed{\phantom{2222}}_{\, 3}$}
\put(3,0){$122\bar{1}$}
\end{picture} \otimes 
\begin{picture}(20,1)
\put(0,0){$\boxed{\phantom{22}}_{\, 6}$}
\put(3,0){$2\bar{2}$}
\end{picture} \otimes
\begin{picture}(30,1)
\put(0,0){$\boxed{\phantom{2222}}_{\, 6}$}
\put(3,0){$22\bar{2}\bar{2}$}
\end{picture} \otimes
\begin{picture}(30,1)
\put(0,0){$\boxed{\phantom{2222}}_{\, 7}\,,$}
\put(3,0){$12\bar{2}\bar{1}$}
\end{picture} 
\quad\quad
\begin{picture}(30,1)
\put(0,0){$\boxed{\phantom{2222}}_{\, 3}$}
\put(3,0){$122\bar{1}$}
\end{picture} \otimes 
\begin{picture}(30,1)
\put(0,0){$\boxed{\phantom{2222}}_{\, 6}$}
\put(3,0){$22\bar{2}\bar{2}$}
\end{picture}\otimes
\begin{picture}(20,1)
\put(0,0){$\boxed{\phantom{22}}_{\, 6}$}
\put(3,0){$2\bar{2}$}
\end{picture} \otimes
\begin{picture}(30,1)
\put(0,0){$\boxed{\phantom{2222}}_{\, 7}\,.$}
\put(3,0){$12\bar{2}\bar{1}$}
\end{picture}
\end{equation*}
Both of them lead to 
\begin{equation*}
p=\Phi_1{\mathcal C}_1(p^{(1)}) = 
\boxed{11}\otimes 
\begin{picture}(18,1)
\put(0,0){$\boxed{\phantom{22}}$}
\put(3,0){$0\bar{0}$}
\end{picture}\otimes
\boxed{22} \otimes \boxed{11} \otimes 
\begin{picture}(18,1)
\put(0,0){$\boxed{\phantom{22}}$}
\put(3,0){$2\bar{2}$}
\end{picture} \otimes \boxed{11} \otimes 
\begin{picture}(18,1)
\put(0,0){$\boxed{\phantom{22}}$}
\put(3,0){$\bar{2}\bar{2}$}
\end{picture} \otimes \boxed{22} \otimes 
\begin{picture}(18,1)
\put(0,0){$\boxed{\phantom{22}}$}
\put(3,0){$0\bar{0}$}
\end{picture} \otimes 
\begin{picture}(18,1)
\put(0,0){$\boxed{\phantom{22}}\,.$}
\put(3,0){$2\bar{2}$}
\end{picture}
\end{equation*}
Thus we obtain 
$\iota^{-1}(p) = 1\,\emptyset\, 2 \, 1 \, 0 \, 1 \, \bar{2} \, 
2 \, \emptyset \, 0, $
in agreement with (\ref{eq:pD(2)}).

\appendix
\section{Crystals and combinatorial R}\label{app:crystal}
The crystals $B_l$ used in the main text
are crystal bases of irreducible finite-dimensional representations of a 
quantum affine algebra $U'_q(\geh)$. 
Let us recall basic facts on them following \cite{Ka,KMN,KKM}.

Let $P$ be the weight lattice, $\{\alpha_i\}_{0\le i\le n}$ the simple 
roots, and $\{\Lambda_i\}_{0\le i\le n}$ 
the fundamental weights of $\geh$.
A crystal $B$ is a finite set with weight decomposition 
$B=\sqcup_{\lambda\in P}B_\lambda$. 
The Kashiwara operators $\tilde{e}_i, \tilde{f}_i$ ($i=0,1,\cdots,n$) 
act on $B$ as 
$
\tilde{e}_i: B_\lambda\longrightarrow B_{\lambda+\alpha_i} \sqcup \{0\},\;
\tilde{f}_i: B_\lambda\longrightarrow B_{\lambda-\alpha_i} \sqcup \{0\}.
$
In particular, these operators are nilpotent.
By definition, we have $\tilde{f}_ib=b'$ if and only if $b=\tilde{e}_ib'$.
For any $b\in B$, set 
$\varepsilon_i(b)=\max\{m\ge0\mid \tilde{e}_i^m b\ne0\}$
and $\varphi_i(b)=\max\{m\ge0\mid \tilde{f}_i^m b\ne0\}$.
Then we have the weight ${\rm wt} b$ of $b$ by 
${\rm wt} b=\sum_{i=0}^n(\varphi_i(b)-
\varepsilon_i(b))\Lambda_i$.

For two crystals $B$ and $B'$, one can define the tensor product
$B\ot B'=\{b\ot b'\mid b\in B,b'\in B'\}$. The operators 
$\tilde{e}_i,\tilde{f}_i$ act on $B\ot B'$ by
\begin{eqnarray*}
\tilde{e}_i(b\ot b')&=&\left\{
\begin{array}{ll}
\tilde{e}_i b\ot b'&\mbox{ if }\varphi_i(b)\ge\varepsilon_i(b')\\
b\ot \tilde{e}_i b'&\mbox{ if }\varphi_i(b) < \varepsilon_i(b'),
\end{array}\right. \\
\tilde{f}_i(b\ot b')&=&\left\{
\begin{array}{ll}
\tilde{f}_i b\ot b'&\mbox{ if }\varphi_i(b) > \varepsilon_i(b')\\
b\ot \tilde{f}_i b'&\mbox{ if }\varphi_i(b)\le\varepsilon_i(b').
\end{array}\right. 
\end{eqnarray*}
Here $0\ot b'$ and $b\ot 0$ should be understood as $0$. 
For crystals we
are considering, there exists a unique isomorphism $B\ot B'
\stackrel{\sim}{\rightarrow}B'\ot B$, {\em i.e.} a unique map
which commutes with the action of Kashiwara operators. In particular,
it preserves the weight.

For a crystal $B$ we define its affinization 
$\Aff(B)=\{b[d]\mid d\in\Z, b\in B\}$ by 
$\tilde{e}_i(b[d])=(\tilde{e}_ib)[d-\delta_{i0}]$ and 
$\tilde{f}_i(b[d])=(\tilde{f}_ib)[d+\delta_{i0}]$. 
($b[d]$ here corresponds to  $T^{-d}af(b)$ in \cite{KMN}.)
The crystal isomorphism 
$B\ot B'\stackrel{\sim}{\rightarrow}B'\ot B$ is lifted up to a map 
$\Aff(B)\ot\Aff(B')\stackrel{\sim}{\rightarrow}\Aff(B')\ot\Aff(B)$
called the combinatorial $R$. It has the following form:
\begin{eqnarray*}
R\;:\;\Aff(B)\ot\Aff(B')&\longrightarrow&
\quad \;\;\;\Aff(B')\ot\Aff(B)\\
b[d]\ot b'[d']\quad\,&\longmapsto&
\tilde{b}'[d'\!-\!H(b\ot b')]\ot \tilde{b}[d\!+\!H(b\ot b')],
\end{eqnarray*}
where $b\ot b'\mapsto\tilde{b}'\ot\tilde{b}$ 
under the isomorphism $B\ot B'
\stackrel{\sim}{\rightarrow}B'\ot B$. $H(b\ot b')$ is called the 
energy function and determined up to an additive constant by
\[
H(\tilde{e}_i(b\ot b'))=\left\{%
\begin{array}{ll}
H(b\ot b')+1&\mbox{ if }i=0,\ \varphi_0(b)\geq\varepsilon_0(b'),\ 
\varphi_0(\tilde{b}')\geq\varepsilon_0(\tilde{b}),\\
H(b\ot b')-1&\mbox{ if }i=0,\ \varphi_0(b)<\varepsilon_0(b'),\ 
\varphi_0(\tilde{b}')<\varepsilon_0(\tilde{b}),\\
H(b\ot b')&\mbox{ otherwise}.
\end{array}\right.
\]
\begin{proposition}[Yang-Baxter equation] \label{prop:YBeq}
The following equation holds on 
$\Aff(B)\ot\Aff(B')\ot\Aff(B'')$:
\[
(R\ot1)(1\ot R)(R\ot1)=(1\ot R)(R\ot1)(1\ot R).
\]
\end{proposition}
We often write the map $R$ simply by $\simeq$.
The combinatorial $R$ is naturally restricted to $B \ot B'$. 

In the main text we are concerned about the crystal 
$B_l$ corresponding to the $l$-fold symmetric fusion 
of the vector representation.
We normalize the energy function so that 
\begin{equation}\label{eq:Hmax}
\max\{H(b\otimes c)\mid b \otimes c \in B_l\otimes B_m\}
= \min(l,m).
\end{equation}
Under this convention one has 
\begin{equation}\label{eq:Hmin}
\min\{H(b\otimes c)\mid b \otimes c \in B_l\otimes B_m\}
= \begin{cases}
0  & \geh_n=A^{(1)}_n, C^{(1)}_n, \\
-\min(l,m) & \geh_n \neq A^{(1)}_n, C^{(1)}_n.
\end{cases}
\end{equation}
When $l=m$, the combinatorial $R$ 
becomes the identity map on $B_l\otimes B_l$ 
but still acts non-trivially as
$R(x[d]\otimes y[e])=x[e-H(x\ot y)]\ot y[d+H(x\ot y)]$.

\begin{acknowledgments}
The authors thank Anne Schilling and Mark Shimozono
for kind interest and comments.
A.K. and Y.Y. are supported by Grants-in-Aid for Scientific 
Research JSPS No.15540363 and No.17340047, respectively.
R.S. is grateful to Miki Wadati for warm encouragement during the
study.
\end{acknowledgments}

\vspace{5mm}
\begin{flushleft}
Atsuo Kuniba:\\
Institute of Physics, Graduate School of Arts and Sciences,
University of Tokyo,
Komaba, Tokyo 153-8902, Japan\\
\texttt{atsuo@gokutan.c.u-tokyo.ac.jp}\vspace{3mm}\\
Masato Okado:\\
Department of Mathematical Science,
Graduate School of Engineering Science,
Osaka University, Osaka 560-8531, Japan\\
\texttt{okado@sigmath.es.osaka-u.ac.jp}\vspace{3mm}\\
Reiho Sakamoto:\\
Department of Physics, Graduate School of Science, 
University of Tokyo, Hongo, Tokyo 113-0033, Japan\\
\texttt{reiho@monet.phys.s.u-tokyo.ac.jp}\vspace{3mm}\\
Taichiro Takagi:\\
Department of Applied Physics, National Defense Academy,
Kanagawa 239-8686, Japan\\
\texttt{takagi@nda.ac.jp}\vspace{3mm}\\
Yasuhiko Yamada:\\
Department of Mathematics, Faculty of Science,
Kobe University, Hyogo 657-8501, Japan\\
\texttt{yamaday@math.kobe-u.ac.jp}
\end{flushleft}
\end{document}